\documentclass[final,leqno]{siamltex}
\pdfoutput=1
\usepackage{amsmath}
\usepackage{amssymb}
\usepackage{amsfonts}
\usepackage{bbm}
\usepackage{graphicx}
\usepackage{booktabs}
\usepackage[section]{placeins}

\newcounter{algonum}
\newcommand{\algo}[1]{\refstepcounter{algonum}\label{#1}

\smallskip
{
\noindent \bf Algorithm \ref{#1}
}
}

\newcommand{\refsimple}{\eqref{eqn-simple1}--\eqref{eqn-simple2}}
\newcommand{\refsf}{\eqref{eqn-sf1}--\eqref{eqn-sf2}}
\newcommand{\refsfwith}
{\eqref{eqn-sf1}--\eqref{eqn-sf2}
with $\ep=0.01$, $a=0.02$, and $\sigma^2=0.113$}
\newcommand{\refepone}
{\eqref{eqn-sf1}--\eqref{eqn-sf2}
with $\ep=1$, $a=5$, and $\sigma^2=0.113$}

\newcommand{\bbOne}{\mathbbm{1}}

\newcommand{\R}{\mathbb{R}}

\newcommand{\bbR}{\mathbb{R}}

\newcommand{\inv}{^{-1}}

\newcommand{\cX}{\mathcal{X}}
\newcommand{\cZ}{\mathcal{Z}}
\newcommand{\cP}{\mathcal{P}}

\newcommand{\cL}{\mathcal{L}}

\newcommand{\ep}{\epsilon}
\newcommand{\lam}{\lambda}

\renewcommand{\P}{{\mathcal{P}}}
\newcommand{\K}{{\mathcal{K}}}
\newcommand{\cK}{{\mathcal{K}}}
\newcommand{\A}{{\mathcal{A}}}

\newcommand{\spann}{\operatorname{span}}

\newcommand{\eng}[2]{\ensuremath{#1 \times 10^{#2} }}

\title{A computational method to extract macroscopic variables and their dynamics in multiscale systems%
}

\author{Gary Froyland \thanks{School of Mathematics and Statistics, University of New South Wales, Australia ({\tt G.Froyland@unsw.edu.au}).}
        \and Georg A. Gottwald  \thanks{School of Mathematics and Statistics, University of Sydney, Australia ({\tt georg.gottwald@sydney.edu.au}).} \and Andy Hammerlindl \thanks{School of Mathematics and Statistics, University of Sydney, Australia and School of Mathematics and Statistics, University of New South Wales, Australia ({\tt andy@maths.usyd.edu.au}).}}

\begin{document}

\maketitle

\begin{abstract}
\noindent
This paper introduces coordinate-independent methods for analysing multiscale dynamical systems using
numerical techniques based on the transfer operator and its adjoint.
In particular, we present a method for testing whether an arbitrary dynamical
system exhibits multiscale behaviour and for estimating the time-scale
separation.
For systems with such behaviour, we establish techniques for analysing the
fast dynamics in isolation, extracting slow variables for the
system, and accurately simulating these slow variables at a large time step.
We illustrate our method with numerical examples and show how the reduced slow dynamics faithfully represents statistical features of the full dynamics which are not coordinate dependent.
\end{abstract}

\begin{keywords}
multiscale systems, slow-fast systems, transfer operator, Koopman operator
\end{keywords}

\begin{AMS}
37xx,  
37Mxx, 
65Pxx 
\end{AMS}

\section{Introduction} \label{sec-intro} 
Devising efficient computational methods to simulate complex
systems is of paramount importance to a wide range of scientific fields, including
biomolecular dynamics, material science and climate
science. The sheer dimensionality, often paired with the presence of a rich
hierarchy of temporal scales and several metastable states, makes direct
computational modelling over the whole range of temporal scales intractable.

Scientists therefore seek reduced models for some designated variables which
carry the relevant information. There are two separate scenarios when such a
reduction is possible: scale separation and weak coupling \cite{givon2004}.
Here we concentrate on the large class of time-scale separated systems which
include models from molecular dynamics, chemical kinetics, and climate. In
scale-separated systems a reduction to slow degrees of freedom has two
computational advantages.
In addition to the obvious computational advantages of decreasing the
dimension, model reductions offer the advantage that the time-step
to be used in simulations can be orders of magnitudes larger since the reduced
model only involves slow variables.

How to extract the set of relevant variables and the
associated relevant dynamics from a dynamical system is one of the most
challenging problems in computational modelling. In this work we address this problem on several levels: We first devise a method which tests whether the complex system under consideration is actually a scale-separated system which allows for a decomposition into slow and fast degrees of freedom. In a second step we then identify possible slow variables. We accurately describe the statistics of the fast dynamics by using transfer operator techniques and then construct the reduced dynamics for the slow variables.

For model reduction techniques to be successful,
a time-scale separation must exist in the first place.
Many treatments simply assume that such behaviour exists.
They further assume the directions of the slow and fast dynamics
are known in advance and that these directions align with the coordinates in
which the system is defined.

In order to analyse general dynamical systems without making such assumptions,
we develop methods based on
the transfer operator and its dual, the so-called Koopman operator.
These two operators
represent a global description of the dynamical
system's action on ensembles and observables, respectively,
and provide powerful tools for
identifying global slow dynamical modes also known as strange eigenmodes,
persistent patterns, metastable sets and almost-invariant sets
\cite{DellnitzJunge99,DeuflhardEtAl00,DeuflhardSchuette04,LiuHaller04,MezicBanaszuk04,mezic2005,Mezic13,Froyland05}.
Transfer operator methodologies have been used to detect
and approximate slowly decaying modes in a number of settings, including
molecular dynamics, ocean dynamics, atmospheric dynamics, and general fluid
flow
\cite{SchuetteEtAl01,FroylandEtAl07,FroylandEtAl10}.

The spectrum of the transfer operator is closely related to the ergodic
properties of the system, the rate of the decay of correlations,
and the overall ``speed'' of the dynamics
\cite{froy1997comp, bkl2002, froy2007}.
To exploit these properties, we develop techniques for isolating the fast
dynamics of the system, computing its transfer operator,
and using this to estimate the scale of the temporal separation.
This study of the fast dynamics in isolation is a topic which has received
comparatively little attention.

Once a clear time-scale separation has been established, one needs to identify the slow variables. We use the eigenfunctions of the Koopman operator for the full system to
define a projection to a lower dimensional space.
This, in effect,
identifies the slow variables of the system.
Our approach is similar to the diffusion map approach \cite{CoifmanLafon06,CoifmanEtAl08,GaveauSchulman06}. Therein slow variables are defined as functions of the physical variables, and the reduction then consists of determining the temporal evolution of those functions.
In contrast, we propose a method that allows for a direct representation of the slow variables, rather than describing the dynamics in function space.

The next step is to construct reduced dynamics for slow variables which reliably approximate the slow dynamics of the full system, with the computational advantages mentioned above. There exists a plethora of methods to construct reduced dynamics for slow variables \cite{givon2004}.
Averaging of deterministic and stochastic systems
over the measure induced by the fast process has been widely used
\cite{VladimirArnold,VerhulstSanders}. In the case when the averaged dynamics
turns out to be trivial, diffusive effects become important and singular
stochastic perturbation theory (homogenization) can be employed
\cite{Khasminsky66,Kurtz73,Papanicolaou76,MajdaEtAl01,JustEtAl01,PavliotisStuart}
to derive reduced equations for the slow variables, even if the underlying
dynamics is deterministic (but sufficiently chaotic)
\cite{MelbourneStuart11,gottwald2013homogenization}. On short time scales,
deterministic and stochastic center manifold theory
\cite{Carr,Boxler89,NamachchivayaEtAl90,NamachchivayaEtAl91,BerglundGentz,Roberts08}
is a well-known method to describe the dynamics locally close to a fixed
point, in which the fast variables are slaved to the dynamics of the slow
variables. A different approach
\cite{ChorinEtAl00,ChorinHald,Stinis07,JustEtAl01} uses the Mori-Zwanzig
projection formalism \cite{MoriEtAl74,Zwanzig60}. 
Some of these methods have a rigorous mathematical footing with clearly stated
assumptions for their validity
\cite{Khasminsky66,Kurtz73,Papanicolaou76,Carr,MelbourneStuart11,gottwald2013homogenization,MacleanGottwald14}.
However, from a practical view point these methods are limited because they require the underlying equations to be of a simple enough form to allow for the required analytical manipulations.

In more complicated situations, i.e., in realistic applications, one has to
resort to numerical simulations.
We mention here the {equation-free} projection algorithms
\cite{GearKevrekidis03,KevrekidisGearEtAl03} and the {heterogeneous multiscale
method} \cite{E03,EEtAl07}.
These methods have been applied to a wide range of
problems, including modelling of water in nanotubes, micelle formation,
chemical kinetics and climate modelling and data
assimilation \cite{MajdaEtAl01,KevrekidisSamaey09,MitchellGottwald12b,GottwaldHarlim13}.
The underlying assumption of these methods is that the fast dynamics quickly relaxes to its equilibrium value (conditioned on the slow variables). The vector field of the reduced slow equation is then estimated by averaging over short bursts of the full dynamics.
The presence of metastable states and the rare transitions between them
severely impede the relaxation towards equilibrium \cite{WalterSchuette06},
rendering short bursts ineffective to estimate the slow dynamics.
Instead we employ the transfer operator to estimate the averages.
The transfer operator has been used to devise reduced slow dynamics
in \cite{Crommelin06,Crommelin11} for estimating
stochastic models from discrete time series. The analytical results in
\cite{Crommelin11} provide support for the approach we choose here.
Our reduction method differs from previous methods
in that we use the invariant measure
computed by the transfer operator
to populate the fast fibers conditioned on the slow states. This cloud of initial conditions is then propagated for an intermediate time scale and, by using the eigenfunctions of the Koopman operator, projected back onto the slow variables. This allows us to estimate the drift and diffusion coefficients via ensemble averages. The existence of the intermediate time scale is assured by homogenisation theory and numerically estimated by requiring that the increments of the slow dynamics are nearly Gaussian.

\medskip

The paper is organised as follows. In Section~\ref{sec-sfabs} we present an abstract framework for slow-fast systems. The basic concepts of the transfer and the Koopman operators are introduced in Section~\ref{sec-ops}. Section~\ref{sec-finding} presents a novel algorithm to effectively detect the presence of a multiscale structure in a given dynamical system as well as identifying the possible slow degrees of freedom. The algorithm is illustrated with a suite of test problems. In Section~\ref{sec-reduced} we present an algorithm to extract the reduced dynamics of the slow variables which were determined in Section~\ref{sec-finding}. Here we use results from homogenization theory to fit parameters of a reduced stochastic differential equation (SDE) for the slow variables. Details on the computational issues involved in calculating the transfer operator and the Koopman operator are discussed in Section~\ref{sec-computing}. We conclude with a discussion and outlook in Section~\ref{sec-discussion}.

\section{An abstract framework for multiscale systems} \label{sec-sfabs} 

In model reduction it is common to consider a projection $\P:\cZ \to \cX$,
where $\cZ$ represents the phase space of the full system, and $\cX$ is a
lower-dimensional space where a reduced model of the dynamics captures its
slowest scale behaviour; see for example \cite{Zwanzig,givon2004}. 
The lower-dimensional space $\cX$ does not necessarily have to be a subspace of $\cZ$. In practice, however, one may want to reduce onto a slow space $\cX \subset\cal{Z}$ to identify slow coordinates as a subset of the original (full) coordinate system. 
In this
paper, we assume that we have no {\em a priori} knowledge about the map
$\P:\cZ \to \cX$, only that it exists.  In particular, we do not assume that
we know what variables are slow or fast, respectively. Therefore, none of the
techniques developed in this paper explicitly use the map $\P$, but rather
exploit the structure it gives to the overall system.

For this, we assume that there is a dynamical system $T:\cZ \to \cZ$ on the full
space and a ``reduced dynamics'' $S:\cX \to \cX$ such that the following
properties hold{:}
\begin{enumerate}
    \item
    $S$ models the long-term behaviour of $T$; that is, there is a large
    iterate $N  \gg  1$ such that $\P \circ T^N \approx S \circ \P$.
    \item
    The short term behaviour of $T$ can be approximated by dynamics on
    a fiber of $\P$.  That is, if
    \[
        z, T(z), \ldots, T^{n}(z)
    \]
    is an orbit of finite length $n  \ll  N$, then $\P(T^i(z))$ is nearly
    constant for $i=0,\ldots,n$ and there is a dynamical system $\hat T:F \to
    F$ on the fiber $F = \P \inv(\P(z))$ and the orbit
    \[
        z, \hat T(z), \ldots, \hat T^{n}(z)
    \]
    of length $n$ closely approximates that of $T$.

\end{enumerate}
\begin{figure}[tp]
\begin{centering}
\includegraphics{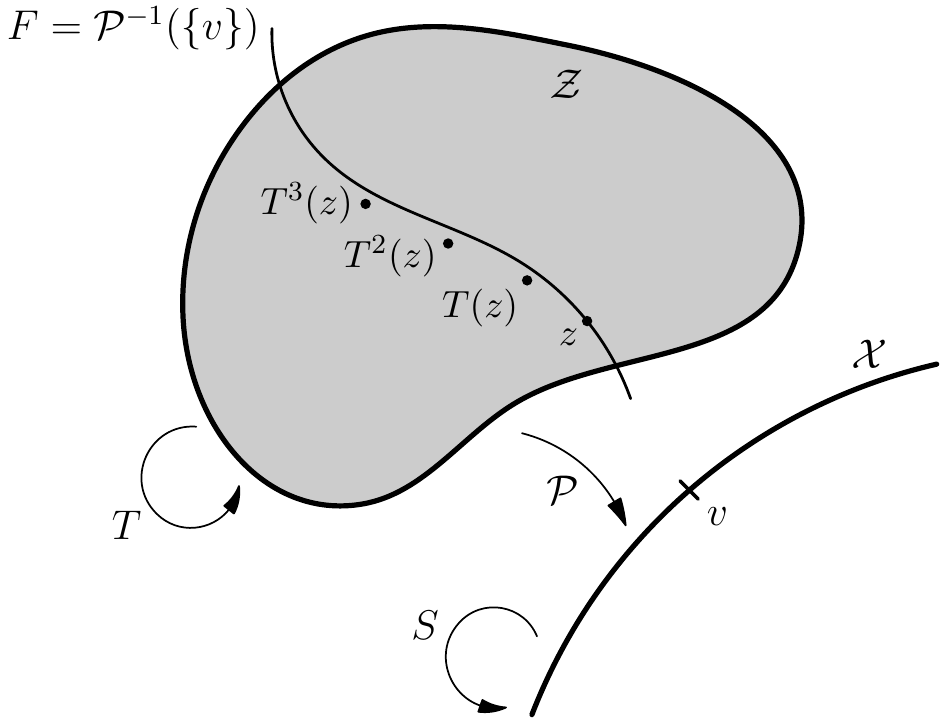}
\end{centering}
\caption{An illustration of the abstract slow-fast framework.
On short time scales, orbits of $T$ stay close to fibers of $\P$.
On long time scales, the orbits of $T$, after projecting by $\P$,
can be approximated by orbits of $S$.}
\label{fig-Proj}
\end{figure}
Figure \ref{fig-Proj} gives an illustration of this abstraction.
Note that there is no unique choice for ${\cal{P}}$.  Indeed, if $h:\cX \to \cX$ is any continuous invertible function, then $h\circ {\cal{P}}$ would also define a map from ${\cal{Z}}$  to ${\cal{X}}$ which satisfies the condition in the slow-fast framework.  However, from a numerical standpoint, one projection may be better than another, and a change of projection is used for such a purpose in Section~\ref{sec-reduced}.\\
For ease of exposition we restricted the presentation here to the case of deterministic systems with discrete time. In the general case, the full dynamics and/or the reduced dynamics could be continuous-time systems, and could have stochastic components, with the obvious modifications.

To turn this abstraction into a fully rigorous definition,
one could consider both the dynamical system
$T_\ep$
and the number $N_\ep$ as depending on a parameter $\ep>0$, such that,
as $\ep$ tends to zero,
\[
    N_\ep \to \infty, \quad \P \circ T_\ep^{N_\ep} \to S \circ \P,
    \quad\text{and}\quad
    T_\ep \to \hat T
\]
under the appropriate notions of convergence.
Indeed, such an approach is used in homogenization theory
\cite{givon2004}.
In this paper, we study numerical examples where $\ep > 0$ is
fixed, and so do not consider these convergence results in detail.
However, for the techniques we develop, a clear time-scale separation is
needed.
The existence of such a time-scale separation is explicitly tested for by our
algorithm.

%
%


\section{The transfer and Koopman operator} \label{sec-ops} 
Information of the long-term behaviour of a dynamical system can often be
obtained from spectral properties of linear operators associated with the
dynamics.
If $T:\cZ \to \cZ$ represents a deterministic, discrete dynamical system,
the transfer operator (or Perron-Frobenius operator) $\cL$ on $L^2(\cZ)$ is
defined by
\begin{equation} \label{eqn-transfer}
    \int g \cdot \cL f \,dm  = \int (g \circ T) \cdot f \,dm
\end{equation}
for all $f,g \in L^2(\cZ)$.  Unless otherwise noted, we take $m$ to be Lebesgue
measure on $\cZ \subset \bbR^d$.
If $\mu$ is a measure defined by a density $f$ (that is $d \mu = f dm$),
then $\cL f $ gives the resulting density after pushing the measure forward by
$T$.  A density $f$ is invariant exactly when $\cL f = f$.
In other words, $f$ is an eigenfunction with eigenvalue $\lam=1$.
If this eigenvalue is simple, then $f$ is the unique invariant density.

In certain settings
\cite{by1993, bkl2002, bosefroydraft},
the transfer operator can be decomposed as $\cL = F + V$
where $F$ is a finite operator and $V$ is contractive, that is, $\|V\| < 1$.
The spectrum of $\cL$ then contains a finite number of isolated eigenvalues
$\lam_i$ satisfying $|\lam_i| > \|V\|$.
In particular, there is a gap in the spectrum between $\lam_1 = 1$ and
$\lam_2$, the next largest eigenvalue in magnitude.
The magnitude of $\lam_2$ can be interpreted as the geometric rate at which
non-equilibrium densities converge to the invariant density,
often called the {\em rate of mixing}.
In general, the leading eigenvalues $\lam_i$ with $|\lam_i|$ close to one
correspond to the slowest decaying transients of the system.
Rigorous results for computing this ``outer spectrum''
of eigenvalues can be found in
\cite{froy1997comp, bkl2002, froy2007}.

A similar theory holds for the adjoint of the transfer operator \cite{MezicBanaszuk04,mezic2005,Mezic13}.  This is
the Koopman operator $\K$ defined simply by $\K g  = g \circ T$.

In practice in this paper, we apply numerical schemes,
which approximate $\cL$ and $\cK$ by finite dimensional operators.
These approximations can be viewed as the
transfer and Koopman operators
$\tilde\cL$ and $\tilde\cK$
now defined for a non-deterministic dynamical system $\tilde T$
which adds local noise to the dynamics of $T$.
The effect of this noise is to ensure that the images of the operators
$\tilde\cL$ and $\tilde\cK$
are finite-dimensional, and so their eigenfunctions exist and can be computed.
These operators can be defined in such a way that the dynamics of $\tilde T$
closely resembles that of $T$, and the numerical implementation
of these operators is described in detail in Section \ref{sec-computing}.
Because of this,
for the remainder of the paper we treat $\cL$ and $\cK$
as if they are compact operators on $L^2(\cZ)$ and assume their leading
eigenvalues and associated eigenfunctions can be calculated numerically.


\section{Testing for multiscale behaviour and identifying fast and slow coordinates} \label{sec-finding} 

We now consider the Koopman operator for a slow-fast system,
and argue that its eigenfunctions can be used to
find the projection $\P:\cZ \to \cX$ without any {\em a priori}
knowledge of the system.

\subsection{Defining fiber dynamics} 
\label{subsec-level}

In the framework for slow-fast systems of Section \ref{sec-sfabs},
we can consider Koopman operators for both
the full dynamics $T$ with
\begin{align*}
\K_Tg=g\circ T\, ,
\end{align*}
and the reduced dynamics $S$ with
\begin{align*}
\K_Sg=g\circ S\, .
\end{align*}
We now compare the spectrum of $\K_S$ with that of $\K_T$.

Suppose that
$\psi:\cX \to \bbR$ is an eigenfunction for the Koopman operator of the reduced dynamics $\K_S$.
That is, $\psi \circ S = \hat \lam \psi$ with eigenvalue $\hat \lam$.
Then,
\[
    \psi \circ \P \circ T^N \approx \psi \circ S \circ \P =
    \hat \lam \,(\psi \circ \P),
\]
so that $\psi \circ \P:\cZ\to\bbR$ approximately solves the eigenfunction
equation $\phi\circ T^N=\hat \lam\phi$ for the Koopman operator of $T^N$.
This suggests that $\K_T$ has an eigenfunction $\phi$ close to $\psi\circ\P$
with eigenvalue $\hat \lam^{\frac{1}{N}}$.
As $N \gg 1$, this eigenvalue is very close to one.

On the other hand, if $\phi:\cZ \to \bbR$ is an eigenfunction of $\K_T$
with eigenvalue $\lam$ such that $0 < |\lam| < 1$,
then as $N$ is large, either $\lam^N$ is close to zero,
or $\lam$ is very close to 1.
In the latter case, for a finite orbit $z, T(z), \ldots T^n(z)$ with $n  \ll  N$
the values $\phi \circ T^j(z) = \lam^j \phi(z)$
are nearly constant.
Since this orbit stays close to the fiber $\P \inv(\{v\})$ and $N$ is large, it
suggests that $\phi$ will be nearly constant along the fiber.
In certain settings, the above arguments can be made rigorous with convergence
results as the time-scale separation tends to infinity
\cite{Crommelin11}.

Suppose the Koopman operator associated with the full dynamics $K_T$ has a number of associated eigenfunctions $\phi_i:\cZ \to \R$
($i=1,\ldots,m$), all with eigenvalues close to one.
Then these have associated
approximations
$\psi_i \circ \P:\cX \to \bbR$.
If, further, $\cX$ can be identified
with its image under the product map
\[
    \psi_1 \times \psi_2 \times \cdots \times \psi_m : \cX \to \bbR^m
\]
then, the product
\[
    \phi_1 \times \phi_2 \times \cdots \times \phi_m : \cZ \to \bbR^m
\]
gives an approximation of the mapping $\P:\cZ \to \cX$.

For simplicity, we restrict our study here to examples where the reduced space $\cX$ is
one-dimensional. In this case the mapping $\P:\cZ \to \cX$ can be well approximated
by a single eigenfunction $\phi:\cZ \to \bbR$.
Possible extensions to higher-dimensional slow subspaces are discussed in Section~\ref{sec-discussion}.

To see if an eigenfunction $\phi$ is actually approximating the (non-linear)
projection $\P$ of a multiscale system, we develop a test for multiscale
behaviour.
%
To do this, we first approximate the fast dynamics $\hat T:F \to F$ defined on
a fiber
$F:=\P \inv(\{v\}) = \phi \inv(\{v\})$.
Assuming $\phi$ is differentiable, its gradient $\nabla \phi$ is a vector
field on $\cZ$.
If this vector field is non-zero on a neighbourhood $U$ of the fast fiber $F$,
it defines a map $\pi$ from $U$ to $F$ simply by flowing
either along $\nabla \phi$ or $-\nabla \phi$ depending on whether the starting
point $z \in U$ satisfies $\phi(z)<v$ or $\phi(z)>v$.
Moreover, one can show by basic calculus
that
\[
    \|\pi(z) - z\|  \le  \frac{1}{C} \|\phi(\pi(z)) - \phi(z)\|
    = \frac{1}{C} \|v - \phi(z)\|
\]
where $C$ is the minimum value of $\|\nabla \phi(z)\|$ on $U$.
For a deterministic $T$, define $\hat T:F \to F$ by
$\hat T = \pi \circ T$.
Then,
\begin{align}
    \|\hat T(z) - T(z)\| < \frac{1}{C} |v - \phi(T(z))|
    = \frac{1}{C} |v - \lambda \phi(z)|
    = \frac{1}{C} |1-\lambda|\ |v|.
\label{eqn-That}
\end{align}

For $\lambda$ close to one, the term $|1-\lam|$ above will be close to zero,
and $\hat T$ will well approximate $T$ and can be used to analyze the fast
dynamics of the system.
The closer to $1$, the better the approximation for the
dynamics $\hat T$ on the fast fibers; therefore, in practice, we use
the eigenfunction $\phi$ which is associated with the second leading eigenvalue
$\lambda_2$
to define the fast fibers and the projection $\P$.

\subsection{Algorithm to test for multiscale behaviour
and identify slow variables} 

To test for multiscale behaviour, we look at the Koopman operator for the
fiber dynamics $\hat T$.  The spectrum of the operator is connected to the
rate of decay of correlations for the system \cite{dellnitz2000}.
Therefore, we expect for the fast dynamics that
the spectrum does not have values as close to one, as for the spectrum
of the full system.

In fact, we propose the following algorithm to test for multiscale behaviour:

\medskip

\algo{algo-compare}
\begin{enumerate}
    \item
    Compute a numerical approximation $\K$ of the Koopman operator of the full
    dynamical system $T$.
    \item
    Determine its leading eigenvalues
    \[
        \lam_1, \lam_2, \lam_3, \ldots
    \]
    such that $1 = |\lam_1|  \ge  |\lam_2|  \ge  |\lam_3|  \ge  \cdots$.
    \item
    Using the eigenfunction $\phi$ associated to $\lam_2$, take a fiber
    $F = \phi \inv(\{v\})$ on which $\|\nabla \phi\|$ is bounded away from zero.
    \item
    Define the ``fiber dynamics'' $\hat T:F \to F$ by $\hat T = \pi \circ T$
    where $\pi$ is determined by the gradient flow associated with $\phi$.
    \item
    Compute a numerical approximation $\hat \K$ of the Koopman operator of the
    fiber dynamics $\hat T$.
    \item
    Determine the leading eigenvalues of $\hat \K$
    \[
        \hat \lam_1, \hat \lam_2, \hat \lam_3, \ldots.  \]
    \item
    Compare $\hat \lam_i$ to $\lam_i$.
    A large ratio implies multiscale behaviour.
\end{enumerate}

\medskip

In a multiscale system, we expect the ratios to be large, roughly on the same
order as the time-scale separation between the slow and fast dynamics, and in
this paper we give examples where this is the case.

If the same techniques are applied to a system without multiscale
behaviour, then the simulated dynamics on a fiber will evolve at
more-or-less the same rate as the full system.
The resulting eigenvalues for the fiber dynamics may be smaller than for
the full system, due to the fiber dynamics acting on a domain of smaller
dimension, but the difference will be modest in nature.
Specific examples of this are given in \ref{sec-simple} and
\ref{sec-epone}.

Once it is determined that a given dynamical system has multiscale behaviour
and that $\phi:\cZ \to \bbR$ reasonably approximates the map $\P:\cZ \to \cX$, we know
that $\hat T$ defined on a fiber reasonably approximates the fast dynamics.
The map $\phi$ in effect defines the slow variable.
It remains to give an approximation for the reduced slow dynamics $S:\cX \to \cX$.
We describe such a construction in Section \ref{sec-reduced} and how it can be
used to simulate the slow dynamics of the system for a time step much larger
than is possible for the full dynamics.

In Algorithm \ref{algo-compare}, one must compute numerical approximations of
Koopman operator and its eigenfunctions.
In the examples in this paper, we use Ulam's method
\cite{ulam-collection,li1976} which approximates the transfer and Koopman
operators by dividing the space into a finite number of boxes.
Care must be taken when applying Ulam's method to multiscale systems, and
specific details about our implementation of the method are given in Section \ref{sec-computing}.
Also, most of our example systems are defined by stochastic differential
equations, and Section \ref{sec-computing} also shows how the operators and
their numerical approximations can be defined in this case.

One could alternatively consider the eigenfunctions of the transfer operator
defined with respect to the invariant measure.
However, as discussed in Section~\ref{sec-computing}, there are significant
numerical problems in computing them to sufficient accuracy.
Therefore, in this paper, we only consider eigenfunctions of the Koopman
operator when computing level sets.

\subsection{Examples} 


\subsubsection{An analytic example} \label{sec-simple} 

Before delving into numerics, we first analyze an overly simple example
which has known analytic formulas for the eigendecomposition.
Consider two independent Ornstein-Uhlenbeck processes given by the SDE
\begin{align}
dx_t &= -x\ dt +\ d{V}_t \label{eqn-simple1} \\
dy_t &= -\frac{1}{\ep^2}y\ dt + \frac{1}{\ep}\,dW_t \label{eqn-simple2}
\end{align}
for some fixed $\ep > 0$ and where $V_t$ and $W_t$ are independent Wiener processes.
If $0<\ep \ll 1$, the $y$-variable evolves much faster than
the $x$-variable and the SDE can be thought of as an example of a multiscale
dynamical system where, in the abstract slow-fast framework of
Section~\ref{sec-sfabs},
$\P:\cZ \to \cX$ is given by
$\P(x,y) = x$.

For now, consider the $x$ variable on its own.  This is a classical
one-dimensional Ornstein-Uhlenbeck system, and its spectral properties are
known \cite{sjogren-survey}.
For each choice of flow time $\tau$ there is an associated Koopman operator
$\K_{\tau}$ for the one-dimensional system given by \eqref{eqn-simple1}.
These form a semi-group of operators $\K_{\tau} = e^{-\tau \A}$, generated by
an operator $\A$.
Note the convention here that we include a minus sign when defining the
generator, implying that eigenvalues of $\A$ are non-negative.
In the appropriate Hilbert space,
the eigenfunctions of $\A$
(and therefore of $K_\tau$ for each $\tau$)
are the Hermite polynomials $H_n$.
The eigenvalues of the infinitesimal generator $\A$ are exactly the
numbers $\{0, 1, 2, \ldots\}$ and $\A H_n = n H_n$ for each $n$.
To keep the analysis independent of the choice of flow-time $\tau$, we express
all eigenvalues in terms of the generator $\A$.

The system in \eqref{eqn-simple2} is the same as that of \eqref{eqn-simple1}, but
with a rescaling of time by a factor of $\ep^{-2}$.  Therefore, its generator
has eigenvalues $\{0, \ep^{-2}, 2\ep^{-2}, \ldots \}$.

For the two-dimensional system, as the two processes are independent,
products of the form
$H_m \times H_n: \bbR^2 \to \bbR,\ (x,y) \mapsto H_m(x) \cdot H_n(y)$
are eigenfunctions, and
the eigenvalues for the full system are then of the form $m + \ep^{-2}n$
for $m,n  \ge  0$.
Assuming $\ep$ is sufficiently small, the leading eigenvalues, i.e., those
with real part closest to zero, will be those with $m$ small and $n$ equal to
zero.
The function $H_0$ is identically equal to one, and therefore
$(H_m \times H_0)(x,y) = H_m(x)$.

The level sets of such a function are lines where the $x$ value is constant,
and
the gradient $\nabla (H_m \times H_0)$ points purely in the $x$ direction.
Therefore, one can verify that the fiber dynamics
$\hat T = \pi \circ T$
described in the last
section corresponds exactly to leaving $x$ constant and evolving $y$ according
to \eqref{eqn-simple2}.  That is, in this simple case, the fiber dynamics
exactly captures the fast dynamics.

The leading eigenvalues associated to the fiber
dynamics are (in terms of the generator) $\{0, \ep^{-2}, 2\ep^{-2}, \ldots \}$
which is exactly $\ep^{-2}$ times the leading eigenvalues of the full system.
Recall that $\ep^{-2}$ is the difference in speed of the slow and fast systems.

In contrast, consider the same system, but with $\ep$ equal to one.
Then, $x$ and $y$ evolve under identical independent processes and one can
check that the eigenspace associated to an eigenvalue $\lam=n$ is of dimension
$n+1$.
In particular, for the second eigenvalue $\lam=1$, every eigenfunction
$\phi: \bbR^2 \to \bbR$ is in the span of $H_1 \times H_0$ and $H_0 \times H_1$
where $H_1(x)=x$ and therefore $\phi$ is
of the form $\phi(x,y)\mapsto (a x+b y)$.
For $\ep = 1$, the system \refsimple\ is invariant under rotation about the
origin, and so we can show that regardless of the exact values of $a$ and $b$,
the fiber dynamics are given by an Ornstein-Uhlenbeck process exactly as
in \eqref{eqn-simple1}. Hence, the eigenvalues associated to the fiber dynamics are
$\{0, 1, 2, \ldots \}$.  The leading ten
eigenvalues of the full system, are (when counted with multiplicity as they
would be if calculated numerically)
\[
    0, 1, 1, 2, 2, 2, 3, 3, 3, 3
\]
and the ratios are
\[
    \frac{0}{0}, \frac{1}{1}, \frac{2}{1}, \frac{3}{2}, \frac{4}{2}, \frac{5}{2}, \frac{6}{3}, \frac{7}{3}, \frac{8}{3}, \frac{9}{3}.
\]
These ratios are growing linearly, but are not large as in the slow-fast case.
Ratios computed for a system with $\ep \approx 1$ will be similar in magnitude.

\subsubsection{A skew-product example} \label{sec-axis} 

We now give an example system which exhibits metastable behaviour in both slow
and fast variables, defined by the following SDE on $\bbR^2$
\begin{align}
dx &= (x - x^3 + \frac{a}{\ep} y)\ dt \label{eqn-sf1} \\
dy &= \frac{1}{\ep^2} (y - y^3)\ dt + \frac{\sigma}{\ep}\,dW_t. \label{eqn-sf2}
\end{align}

Here, $W_t$ represents the standard Wiener process.
Homogenization results (see \cite{givon2004,PavliotisStuart})
show that as $\ep \to 0$, the
evolution of the $x$ variable
approaches that of a drift-diffusion process with the contribution of
$\frac{a}{\ep} y$ replaced by Gaussian noise.
For numerical analysis, however, we only consider the system with the
constants fixed at
$\ep=0.01$, $a=0.02$, and $\sigma^2 = 0.113$.
\begin{figure}[tp]
\begin{centering}
\includegraphics[scale=0.8]{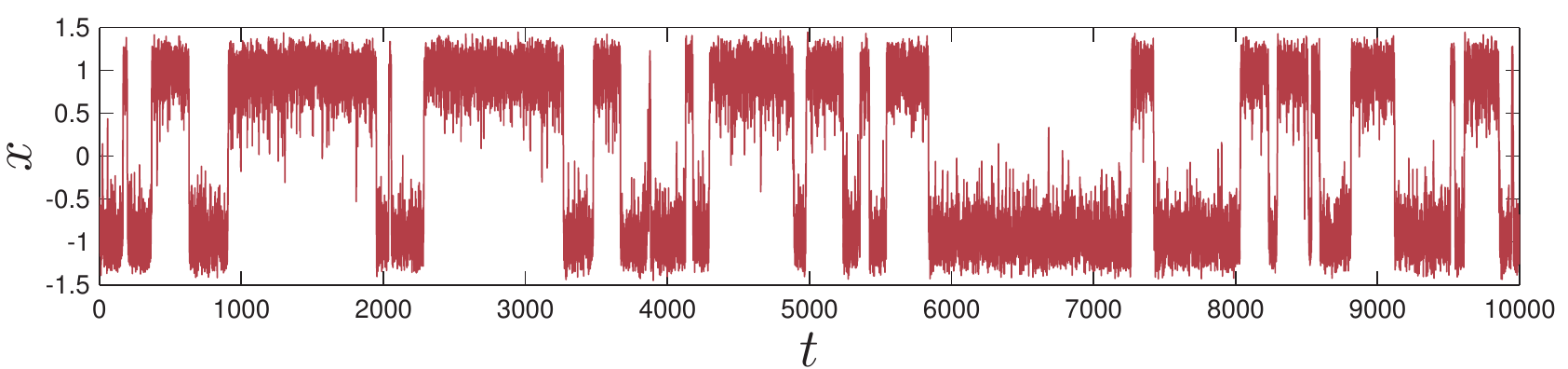}
\includegraphics[scale=0.8]{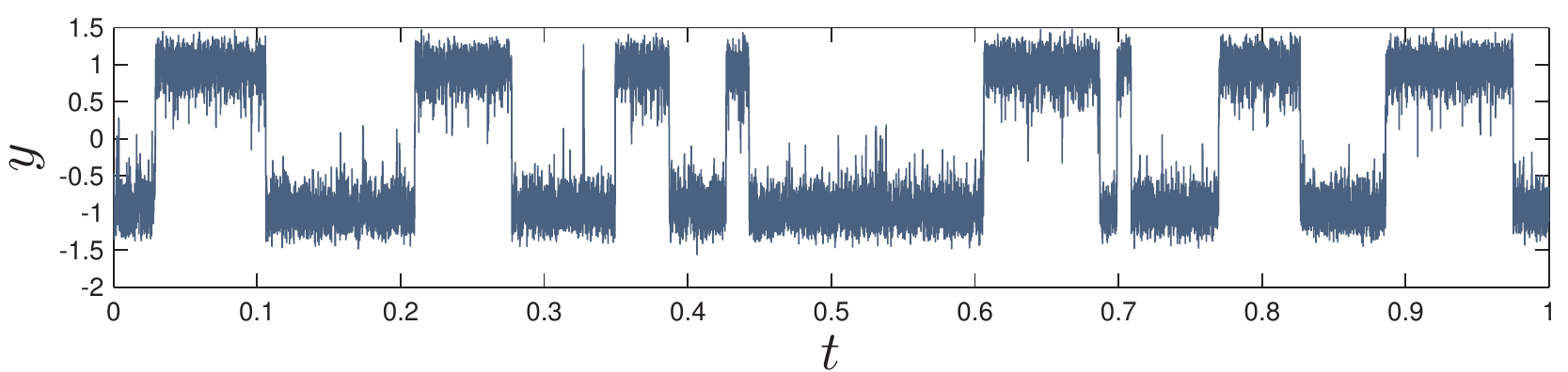}
\end{centering}
\caption{The slow variable $x$ and fast variable $y$ computed for
the system \refsfwith.  The $x$ variable is plotted over an interval
of $10^4$ units of time, while the $y$ variable is plotted over an interval
$1$ unit of time.}
\label{fig-fullevolution}
\end{figure}
For these parameters the fast $y$ variable switches between
states near $y=-1$
and $y=+1$ on the order of about 0.1 units of time.
The slow $x$ variable switches between states near $x=-1$
$x=+1$ on the order of hundreds of units of time, as illustrated
in Figure \ref{fig-fullevolution}.
Because of the long time step needed to switch between them, we refer to these
two states near $x=\pm 1$ as the {\em slow metastable states}
of the system.

For the algorithms developed here, the dynamical system is treated as a
black box.  It is represented by a computer routine which, given initial data
$z = (x,y) \in \bbR^2$, returns a point $T(z) \in \bbR^2$ given by the
Euler-Maruyama integration method
\cite{kloeden1995}
after a single time step of
$\Delta t = \eng{2}{-6}$.
As the system is defined by a stochastic differential equation, calling the
routine multiple times with the same input will yield different outputs.

The value $\Delta t$ was chosen to be large enough to allow efficient
computation of the system while being small enough to still simulate the fast
dynamics of the system to an acceptable degree of accuracy.

\begin{figure}[tp]
\includegraphics{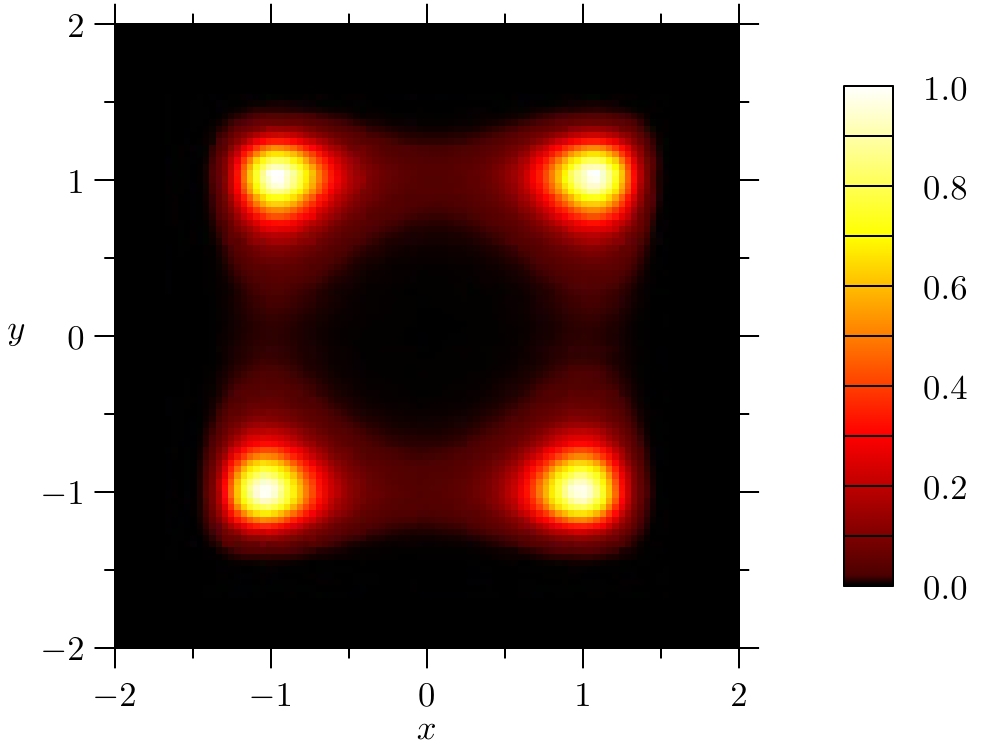}
\caption{The invariant density of the system \refsfwith,
computed using Ulam's method with a flow time of $\tau=0.04$.
(See Section \ref{sec-computing} for numerical details.)}
\label{density-figure}
\end{figure}

\begin{figure}[tp]
\includegraphics{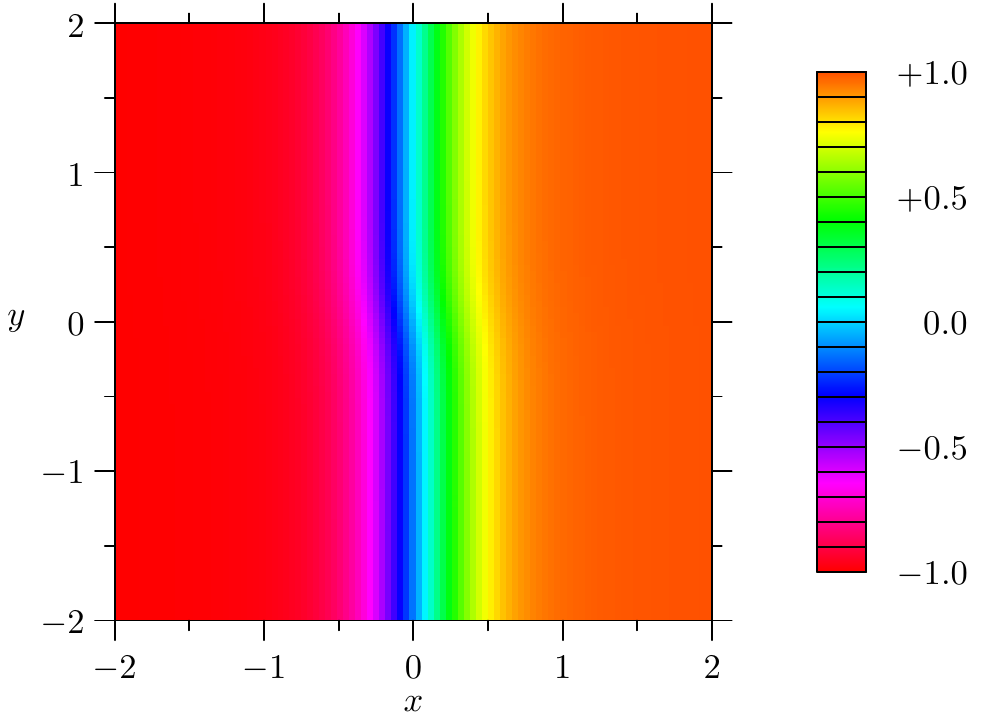}
\caption{The second eigenfunction of the Koopman operator for the system
\refsfwith,
computed using Ulam's method with a flow time of $\tau=0.04$.
(See Section \ref{sec-computing} for numerical details.)}
\label{evec-figure}
\end{figure}

To compute approximations of the eigenfunctions using Ulam's method, we apply
$\eng{2}{4}$ iterates of this time step $\Delta t=\eng{2}{-6}$ for a
total time of $\tau = 0.04$.
Note that while this time interval corresponds to a large number of iterates,
it is still much smaller than the time necessary for $x$ to switch between its
two metastable states.
However, Ulam's method computes accurate approximations of the eigenfunctions
as shown in Figures \ref{density-figure} and \ref{evec-figure}.

We now choose a fiber and simulate dynamics restricted to that fiber
as described in Algorithm \ref{algo-compare}.
As will be explained in Section \ref{sec-computing}, a time interval of
$\hat \tau = \eng{4}{-5}$ is sufficient to accurately
compute eigenvalues of the fiber dynamics.
Note that this $\hat \tau$ is much smaller than the time of $\tau=0.04$ used for
the full dynamics.

\begin{table}[ht]
\caption{Leading eigenvalues $\lam_i$ of the Koopman operator for the full
two-dimension dynamics of the system \refsfwith\ with flow time $\tau=0.04$.
Also computed are the flow-independent values
$\chi_i=-\frac{1}{\tau} \log |\lam_i|$ corresponding to the infinitesimal generator
of the Koopman operator.
Next, are the leading eigenvalues of the Koopman operator of the fiber
dynamics computed with flow time $\hat \tau = \eng{4}{-5}$
and the corresponding $\hat \chi_i=-\frac{1}{\hat \tau} \log |\hat \lam_i|$.
Here, the fiber $\phi \inv(\{v\})$ is for $v=0.8$ where $\phi$ is the
eigenfunction plotted in Figure \ref{evec-figure}.
The final column shows the ratio $\hat \chi_i/\chi_i$.}
\centering
\begin{tabular}{l l l l l l r}
\toprule
\multicolumn{2}{c}{Full dynamics} & \phantom{abc} &
\multicolumn{2}{c}{Fiber dynamics} & \phantom{abc} &
\multicolumn{1}{c}{Ratio} \\
\cmidrule{1-2} \cmidrule{4-5}
\multicolumn{1}{c}{$\lam_i$} &
\multicolumn{1}{c}{$\chi_i$} &&
\multicolumn{1}{c}{$\hat \lam_i$} &
\multicolumn{1}{c}{$\hat \chi_i$} &&
\multicolumn{1}{c}{$\hat \chi_i / \chi_i$}
\\
\hline
1.0000  & 0 && 1.0000 & 0 && \multicolumn{1}{c}{---} \\
0.9995  & \eng{1.25}{-2} && 0.9979 & \eng{5.26}{1} && \eng{4.20}{3} \\
0.9665  & \eng{8.52}{-1} && 0.7212 & \eng{8.17}{3} && \eng{9.59}{3} \\
0.9429  & 1.47 && 0.5838 & \eng{1.35}{4} && \eng{9.15}{3} \\
0.9153  & 2.21 && 0.4694 & \eng{1.89}{4} && \eng{8.55}{3} \\
0.8801  & 3.19 && 0.3463 & \eng{2.65}{4} && \eng{8.30}{3} \\
0.8397  & 4.37 && 0.2424 & \eng{3.54}{4} && \eng{8.11}{3} \\
0.7954  & 5.72 && 0.1619 & \eng{4.55}{4} && \eng{7.95}{3} \\
0.7477  & 7.27 && 0.1055 & \eng{5.62}{4} && \eng{7.74}{3} \\
0.6969  & 9.02 && 0.0646 & \eng{6.85}{4} && \eng{7.59}{3} \\
\bottomrule
\end{tabular}
\label{table-sf}
\end{table}
The leading ten eigenvalues of both the full dynamics and the fiber
dynamics are given in Table \ref{table-sf}.
For this example, a range of several values of $v$ was used to compute dynamics on the
fibers $\phi \inv(\{v\})$ and gave similar eigenvalues.
As these were computed using the different time steps $\tau$ and $\hat \tau$,
we also give the eigenvalues $\chi_i$ for the infinitesimal generator, computed
by taking $\chi_i = -\tfrac{1}{\tau} \log(\lam_i)$.
The final column gives the ratio $\hat \chi_i/\chi_i$ of the eigenvalues of the
two generators.
These ratios are large (over $10^3$), and all within the same order of
magnitude,
which gives strong evidence the system has slow-fast behaviour with a
time-scale separation of between $10^3$ and $10^4$.
Compare this to the scaling term $\ep^{-2} = 10^4$ in equation \eqref{eqn-sf2}
and with the observed switching times for $x$ and $y$ in Figure
\ref{fig-fullevolution}.

As a further test that the fiber dynamics is a valid approximation of the
fast dynamics, we populated the fiber with $10^4$ points $\{z_i\}$
distributed in proportion to the computed invariant density of the system, and
then computed $\|T(z_i) - z_i\|$ and $\|\hat T(z_i) - T(z_i)\|$
for each point. For this example, the average value of $\|T(z_i) - z_i\|$ was $\eng{3.69}{-2}$
and the average value of $\|\hat T(z_i) - T(z_i)\|$ was $\eng{4.69}{-4}$.
Further, the maximum value of $\|\hat T(z_i) - T(z_i)\|$ over all of the points
$z_i$ was $\eng{1.26}{-2}$.  This shows that the fiber dynamics $\hat T$
closely approximates $T$.


\FloatBarrier 

\subsubsection{A ``slow-slow'' system} \label{sec-epone} 

We now show that the same procedure when applied to a system without
multiscale behaviour gives noticeably different results.
As with the last example, consider the system defined by \refsf, but
now set $\ep = 1$ and $a = 5$.  The value $\sigma^2=0.113$ is as before.
Under these parameters, $y$ stays near $+1$ or $-1$ for long periods of
time and, because of the term $\frac{a}{\ep} y$ in \eqref{eqn-sf1},
the $x$ variable is now highly correlated with $y$.
This system does not exhibit multiscale behaviour.

As with the last example, we computed eigenfunctions and then populated a
fiber $\phi\inv(\{v\})$ with $10^4$ points
$\{z_i\}$ distributed in proportion to the computed invariant density of the
system, and then computed $\|T(z_i) - z_i\|$ and $\|\hat T(z_i) - T(z_i)\|$
for each point.
In this case, the average value of $\|T(z_i) - z_i\|$ was $\eng{5.36}{-3}$
and the average value of $\|\hat T(z_i) - T(z_i)\|$ was $\eng{5.34}{-3}$.
These values are nearly equal, showing that the full dynamics does not stay
close to the fiber on short time scales, and that $\hat T = \pi \circ T$ is
not a valid approximation of the dynamics.

Applying Algorithm \ref{algo-compare} produces eigenvalues as shown in
Table \ref{table-epone}.  One should not read too much into the table, as the
fiber dynamics do not meaningfully correspond to any form of ``fast''
dynamics, but notice that the ratios $\hat \chi_i / \chi_i$ are not large in
general and that the first ratio $\hat \chi_2 / \chi_2 \approx 197.37$ is of
significantly different order than the others.

\begin{table}[ht]
\caption{Leading eigenvalues $\lam_i$ of the Koopman operator for the full
two-dimensional dynamics of the ``slow-slow'' system \refepone\ with flow time
$\tau=0.04$.
Also computed are the flow-independent values
$\chi_i=-\frac{1}{\tau} \log |\lam_i|$ corresponding to the infinitesimal generator
of the Koopman operator.
Next, are the leading eigenvalues of the Koopman operator of the fiber
dynamics also computed with flow time $\hat \tau = 0.04$
and the corresponding $\hat \chi_i=-\frac{1}{\hat \tau} \log |\hat \lam_i|$.
The final column shows the ratio $\hat \chi_i/\chi_i$.}
\centering
\begin{tabular}{l l l l l l r}
\toprule
\multicolumn{2}{c}{Full dynamics} & \phantom{ab} &
\multicolumn{2}{c}{Fiber dynamics} & \phantom{ab} &
\multicolumn{1}{c}{Ratio} \\
\cmidrule{1-2} \cmidrule{4-5}
\multicolumn{1}{c}{$\lam_i$} &
\multicolumn{1}{c}{$\chi_i$} &&
\multicolumn{1}{c}{$\hat \lam_i$} &
\multicolumn{1}{c}{$\hat \chi_i$} &&
\multicolumn{1}{c}{$\hat \chi_i / \chi_i$}
\\
\hline
1.0000  & 0 && 1.0000 & 0 && \multicolumn{1}{c}{---} \\
0.9997           & \eng{7.50}{-3} && 0.9425           &  1.48 &&  197.37 \\
0.9674           & \eng{8.29}{-1} && 0.8105           &  5.25 &&    6.34 \\
0.9479           & 1.34           && 0.7019           &  8.85 &&    6.62 \\
0.9470           & 1.36           && 0.5365           &  15.6 &&   11.43 \\
0.9254           & 1.94           && $0.5090 + 0.0237i$ & 16.9  &&    8.70 \\
$0.8983 + 0.0411i$ & 2.66         && $0.5090 - 0.0237i$ & 16.9  &&    6.39 \\
$0.8983 - 0.0411i$ & 2.66         && $0.4903 + 0.0814i$ & 17.5  &&    6.58 \\
0.8961           & 2.74           && $0.4903 - 0.0814i$ & 17.5  &&    6.37 \\
0.8647           & 3.63           && $0.4676 + 0.1403i$ & 17.9  &&    4.93 \\
\bottomrule
\end{tabular}
\label{table-epone}
\end{table}

\subsubsection{A distorted example} \label{sec-distort} 
\label{subsec-dist1}

A major advantage of our approach to identifying and analyzing slow-fast
systems is that the techniques are independent of the choice of coordinates
and therefore the slow and fast directions do not need to align with
coordinate axes in order for the analysis to work.  To demonstrate this, we
apply a transformation to system \refsf.  Recall that in Section
\ref{sec-axis}, we treated the SDE \refsf\ as a black box.
Given a point $z \in \bbR^2$, we have an opaque method for computing the point
$T(z)$ which is the evolution after a time-step
$\Delta t = \eng{2}{-6}$.

\begin{figure}[tp]
\begin{minipage}[b]{0.5\linewidth}
\includegraphics[scale=0.9]{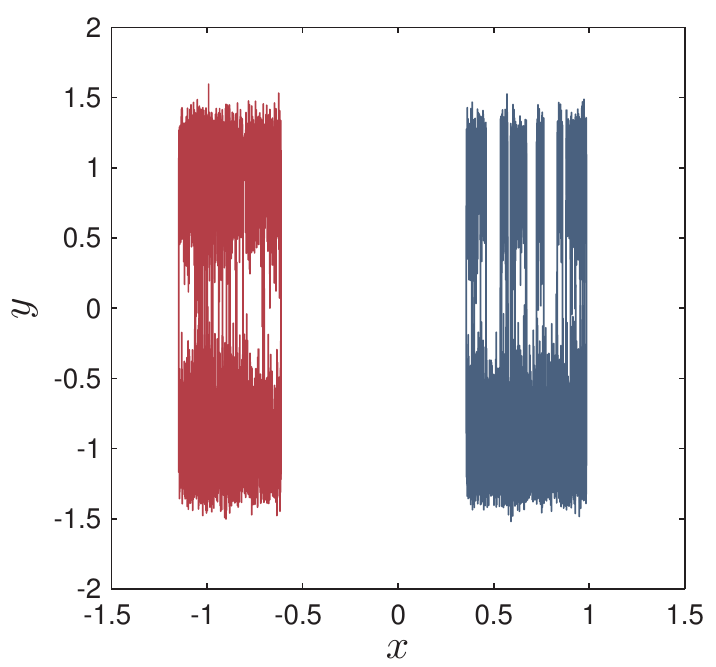}
\end{minipage}
\quad
\begin{minipage}[b]{0.5\linewidth}
\includegraphics[scale=0.9]{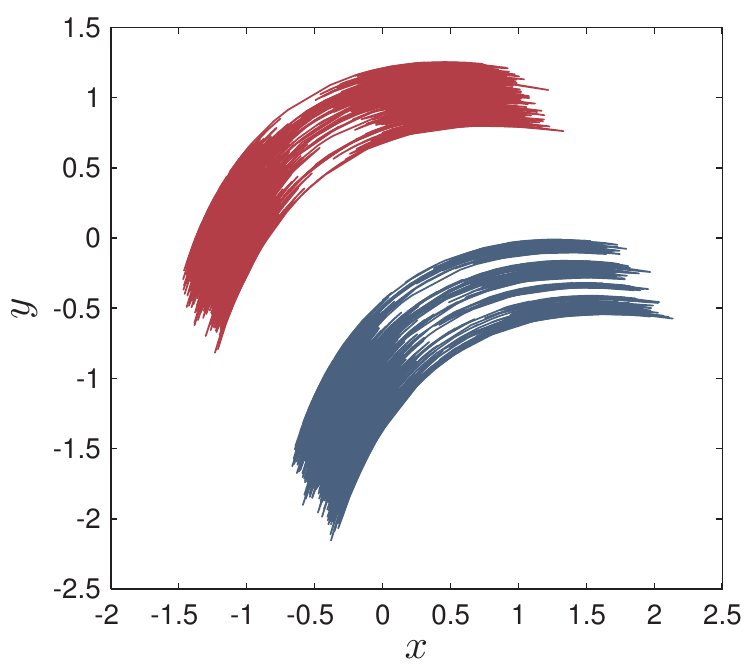}
\end{minipage}
\caption{Two finite orbits of the system \refsfwith\ before and after applying the
distortion $g$ described in Section~\ref{sec-distort}.
Each orbit represents flowing for one unit of time,
and the two orbits were chosen to demonstrate the behaviour of the system near
the two slow metastable states. }
\label{fig-distorb}
\end{figure}

In place of $z \mapsto T(z)$ consider instead,
$z \mapsto g \circ T \circ g \inv (z)$
where $g$ is a diffeomorphism of the domain $\bbR^2$.
This changes the coordinates without actually changing the dynamics.
We take $g$ to be the composition $h \circ R$,
where $h(x,y) = (x + 0.3\, y^2, y)$ and $R$ is an irrational rotation of the
plane (specifically, a clockwise rotation by one radian).
Figure \ref{fig-distorb} shows this distortion applied to two orbits of the system.
\begin{figure}[tp]
\includegraphics{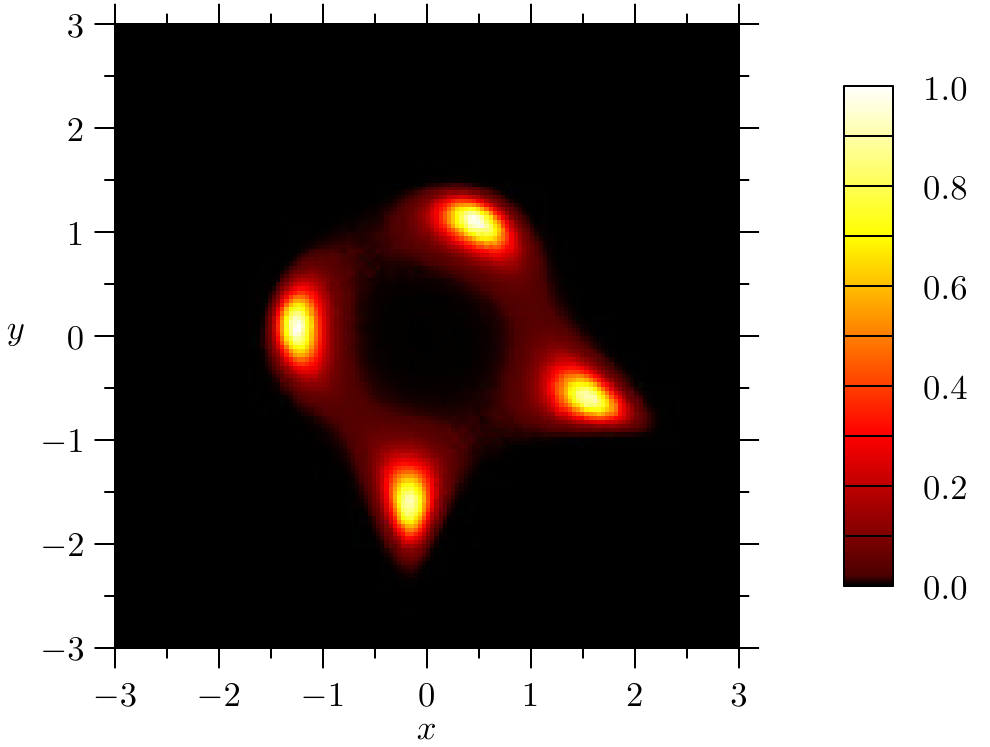}
\caption{The invariant density of the system \refsfwith\ after
the distortion $g$ described in Section \ref{sec-distort}.
The density is computed
by Ulam's method applied to the distorted system using a flow time of
$\tau=0.04$.
}
\label{fig-dist-density}
\end{figure}

\begin{figure}[tp]
\includegraphics{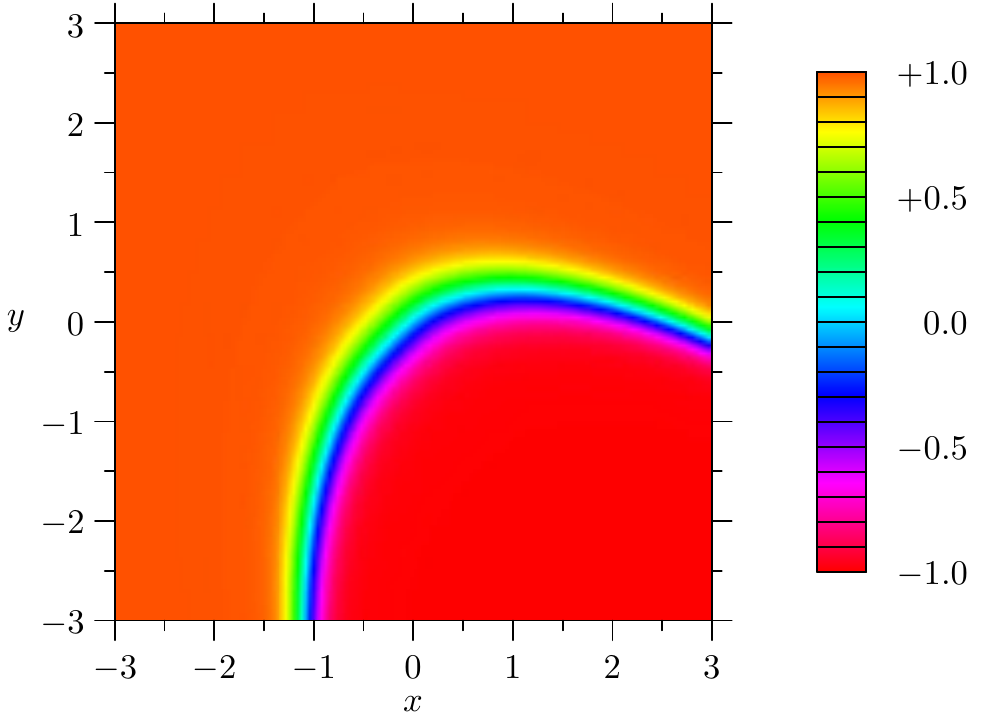}
\caption{The second eigenfunction of the Koopman operator for the system \refsfwith\ after
the distortion $g$ described in Section \ref{sec-distort}.
The eigenfunction is computed
by Ulam's method applied to the distorted system using a flow time of
$\tau=0.04$.
}
\label{fig-dist-evec2}
\end{figure}
It is not immediately clear in
the distorted picture where exactly the slow and fast directions lie.  Still,
we can apply Ulam's method to calculate the eigenfunctions of the transfer
and Koopman operators, as illustrated in Figures \ref{fig-dist-density} and
\ref{fig-dist-evec2}.
The leading eigenvalues should, of course, be the same as in the
original system, and computation confirms that, as shown in
Table \ref{table-dist}.

\begin{table}[]
\caption{Leading eigenvalues for the Koopman operator of the full system
and fiber dynamics as described in Algorithm \ref{algo-compare}.
These eigenvalues are computed using Ulam's method with axis-aligned
boxes for the original system \refsfwith\ and for the same system conjugated
by the diffeomorphism $g:\bbR^2 \to \bbR^2$ defined in Section \ref{sec-distort}.}
\centering
\begin{tabular}{c c c c c}
\toprule
\multicolumn{2}{c}{Full Dynamics} & \phantom{ab} &
\multicolumn{2}{c}{Fiber Dynamics} \\
\multicolumn{2}{c}{($\tau = 0.04$)} &&
\multicolumn{2}{c}{($\hat \tau = \eng{4}{-5}$)} \\
\cmidrule{1-2} \cmidrule{4-5}
original &
distorted &&
original &
distorted
\\
\hline
1.0000  & 1.0000 && 1.0000 & 1.0000 \\
0.9995  & 0.9994 && 0.9979 & 0.9980 \\
0.9665  & 0.9663 && 0.7212 & 0.7230 \\
0.9429  & 0.9423 && 0.5838 & 0.5838 \\
0.9153  & 0.9148 && 0.4694 & 0.4675 \\
0.8801  & 0.8787 && 0.3463 & 0.3461 \\
0.8397  & 0.8385 && 0.2424 & 0.2423 \\
0.7954  & 0.7939 && 0.1619 & 0.1615 \\
0.7477  & 0.7454 && 0.1055 & 0.1027 \\
0.6969  & 0.6944 && 0.0646 & 0.0634 \\
\bottomrule
\end{tabular}
\label{table-dist}
\end{table}
Using the second eigenfunction of the Koopman operator, we compute fiber
dynamics as before.  Eigenvalues associated to the fiber
dynamics are computed and are nearly identical to the eigenvalues computed for
the undistorted system.
This shows that this technique of comparing eigenvalues to
detect multiscale systems is independent of the choice of coordinates.

As before, we tested the validity of the fiber dynamics
by populating a fiber with $10^4$ points $\{z_i\}$
distributed in proportion to the computed invariant density of the system, and
then computed $\|T(z_i) - z_i\|$ and $\|\hat T(z_i) - T(z_i)\|$
for each point.
For the distorted system, the average value of $\|T(z_i) - z_i\|$ was
$\eng{4.37}{-2}$ and the average value of $\|\hat T(z_i) - T(z_i)\|$ was
$\eng{2.99}{-4}$.  Further, the maximum value of $\|\hat T(z_i) - T(z_i)\|$
over all of the points $z_i$ was $\eng{1.26}{-2}$.  This shows that, as in the
undistorted case, the fiber dynamics $\hat T$ closely approximates $T$.


\section{Reduced dynamics} \label{sec-reduced} 

The fiber dynamics approach studied in the previous sections gives a good way to
isolate the fast dynamics of the system.  We also want to efficiently emulate
the slow dynamics of the system using a large time step.  To do this, we
assume the slow dynamics of the system is well approximated by a
drift-diffusion process.
The details and justification of this approach will be given later in the
section, but we first list the proposed algorithm here{:}

\algo{algo-reduce}
\begin{enumerate}
    \item
    Numerically compute the invariant density $\rho$ of the dynamical system
    $T$.
    \item
    Numerically compute an approximate eigenfunction $\phi$ of the Koopman
    operator associated to an eigenvalue close to one.
    \item
    Reparameterize $\phi$ to yield a function $\theta$ which varies evenly
    throughout the domain.
    \item
    Compute fibers for $\theta$.
    \item
    For each computed fiber $\theta \inv(\{v\})$:
    \begin{enumerate}
        \item
        Populate $\theta \inv(\{v\})$ by an ensemble of points
        $\{z_i\}_{i=1}^q$ sampled in proportion to the invariant density
        $\rho$\footnote{Alternatively one could use the fast fiber
        dynamics defined in Section~\ref{subsec-level} to compute the
        invariant density on each fast fiber $F$. However, this is
        computationally much more involved and we therefore use the invariant
        density of the full dynamics.} .
        \item
        Iterate each point $z_i$ forward by $T^k$ and project by $\theta$.
        \item Compute $\alpha_k(v)$ and $\beta_k(v)$
        as the mean and variance of the ensemble
        $\{\vartheta_i\}_{i=1}^q = \{\theta(T^k(z_i)) - v\}_{i=1}^q$.
        That is,
        \[
            \alpha_k(v) = \frac{1}{q}\sum_i \vartheta_i
            \quad\text{and}\quad
            \beta_k(v) = \frac{1}{q}\sum_i
            \bigl(\vartheta_i - \alpha_k(v))^2.
        \]
        \end{enumerate}
    \item
    Extend $\alpha_k$ and $\beta_k$ by linear interpolation to functions on
    $\R$.
    \item
    Use $\alpha = \alpha_k / k \Delta t$ and $\beta = \beta_k / k \Delta t$ to
    define a drift-diffusion process
    \[
        dX_t = \alpha(X_t) dt + \sqrt{\beta(X_t)} dW_t
    \]
    representing the reduced slow dynamics.
\end{enumerate}
\begin{figure}[tp]
\begin{centering}
\includegraphics{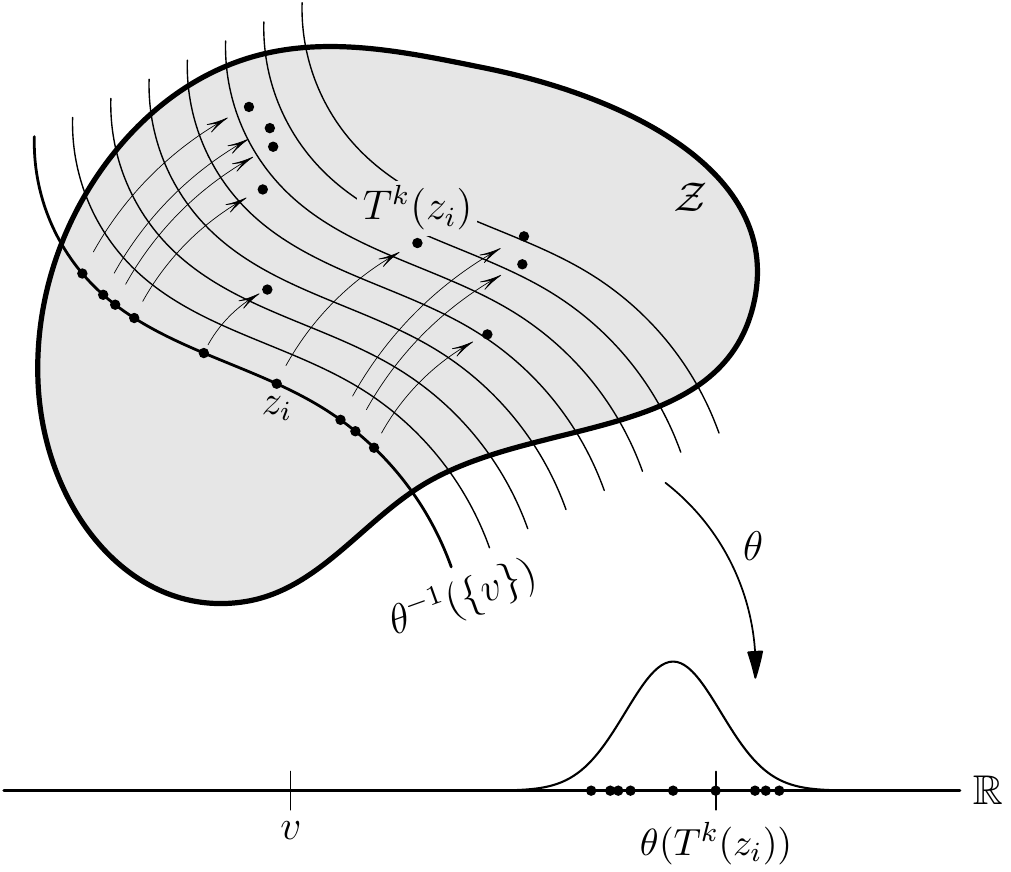}
\end{centering}
\caption{An illustration depicting Algorithm \ref{algo-reduce}.}
\label{fig-Reduce}
\end{figure}
A visual depiction of this algorithm is given in Figure \ref{fig-Reduce}.\\

Let us first comment on Step 3 of the algorithm and how we define the projection $\theta$. Since we are assuming fast dynamics occurs along the level sets of an
eigenfunction $\phi$ of the Koopman operator, this function $\phi$ presents a
natural choice for projection to a one-dimensional space for the slow
dynamics, as discussed in Section~\ref{subsec-level}. However, this is not an ideal candidate for numerical calculation in systems with meta-stable slow states: Suppose the dynamical system has two metastable states in the slow direction,
as in the equations \refsf.
The eigenfunction $\phi$, after scaling, will
be close to $+1$ near one of the metastable states, close to $-1$ near the
other, and change sharply in the area where it transitions from one state to
the other.  Such a sharp change is visible in Figure \ref{evec-figure}.
In fact, such regions where an eigenfunction is nearly constant can be used to
detect the existence of such metastable states and of ``almost invariant
sets'' \cite{dellnitz1999, froypad2009}. For our purpose, however, those regions make it hard numerically to use $\phi$ as a projection.\\
To overcome this problem, we replace $\phi$ by a function which has
similar level sets, but which varies evenly throughout the phase space.\\ Suppose $\phi$ is represented numerically by its values on a finite grid of
points $\{z_i \}_{i=1}^{M}$ in the phase space $\cZ$.
Then, we define the function
$\theta: \{z_i \}_{i=1}^{M} \to \{1, \ldots, M\}$ as a bijection of finite sets and
such that
$\phi(z_i) < \phi(z_j)$ implies
$\theta(z_i) < \theta(z_j)$.
That is, from $\phi$ we have extracted only its ordering of points on the grid.
To get a suitable continuous function,
extend $\theta$ to the rest of the phase space by linear interpolation in
each coordinate.
Let $B \subset \bbR^d$ be the box on which the grid of points lies.
Then, one way of thinking about this step, is that we are replacing $\phi$ by some
reparameterization $\theta = h \circ \phi$ with $h:\bbR \to \bbR$ a
homeomorphism, so that for an interval $[a,b] \subset \bbR$ the
volume of $\theta \inv([a,b]) \cap B$ is roughly proportional to the length of
$[a,b]$.\\
For example, the
eigenfunction plotted in Figure \ref{evec-figure} was computed using a 200 by
200 grid $G$ of points $z_i=(x_i,y_i) \in \bbR^2$.  The relative ordering of the
values $\phi(x_i,y_i)$ on this grid $G$ defines a bijection
$\theta:G \to \{1, \ldots, 40\,000\}$ which extends to the convex hull of $G$ by
bilinear interpolation.  Such a construction would yield a $\theta$ ranging
between $1$ and $40\,000$.  As a further step, we rescale $\theta$ so that it ranges between $-4$ and $+4$.
This makes it easier to allow for an easier comparison between $\theta(T^k(z_i))$ and the original slow $x$-variable in \refsf.\\

Algorithm~\ref{algo-reduce} relies on the reduced system being a stochastic differential equation. Diffusion limits of stochastic multiscale systems such as \refsf\ have been rigorously established in the framework of homogenization (see, for example, \cite{givon2004,PavliotisStuart}) in the limit $\ep\to 0$. Diffusion limits for deterministic dynamical systems have also been established \cite{MelbourneStuart11,gottwald2013homogenization} and have been extended to the discrete time case as well \cite{gottwald2013homogenization}. Assuming,
then, that the slow dynamics is well modelled by a one-dimensional
drift-diffusion system, our algorithm describes how to estimate the drift term
$\alpha(X)$ and the diffusion term $\beta(X)$ in
\begin{equation} \label{eqdd}
    dX_t = \alpha(X_t) dt + \sqrt{\beta(X_t)} dW_t
\end{equation}
so that an orbit of this process closely resembles the projection under $\theta$
of an orbit of the full system.\\ We estimate $\alpha$ and $\beta$ at $v \in \bbR$ by averages over an ensemble
\[
\{X_{t_k,i}\}_{i=1}^{q}
\]
of size $q$. The ensemble is generated as follows (cf. Figure~\ref{fig-Reduce}): Populate the fiber $F$ by a number of points $z_i$, evolve each point forward by $k$ time-steps of the computation to a point to $T^k(z_i)$ at time $t_k$, and then project onto the slow subspace using the projection $\theta$ to yield $\theta(T^k(z_i))$. The initial population of points $\{z_i\}_{i=1}^q$ on $F$ should be chosen in proportion to the invariant density function computed for the full dynamics (or the dynamics on the fast fibers). We remark that this procedure bears resemblance with the equation-free approach proposed in \cite{GearKevrekidis03,KevrekidisGearEtAl03,KevrekidisSamaey09,CoifmanLafon06,CoifmanEtAl08}.\\
%
%
%
%
%
%
The drift and diffusion coefficients $\alpha(v)$ and $\beta(v)$
are then approximated for small times $t_k$ by $\alpha_k(v)/t_k$ and $\beta_k(v)/t_k$
where $\alpha_k(t_k)$ and $\beta_k(t_k)$ are the mean and variance of
$\{X_{t_k,i} - v\}_{i=1}^{q}$ \cite{Gardiner}.\\
%
%
%
%
%
Note that in general the variance of $\{X_{t_k,i} - v\}_{i=1}^{q}$ would
depend on the drift term as well as the diffusion term in \eqref{eqdd};
however
for small times $t_k$ the contribution of the drift term is order $t_k$
and can be neglected with respect to the diffusion term.

\begin{figure}[tp]
\begin{center}
\includegraphics{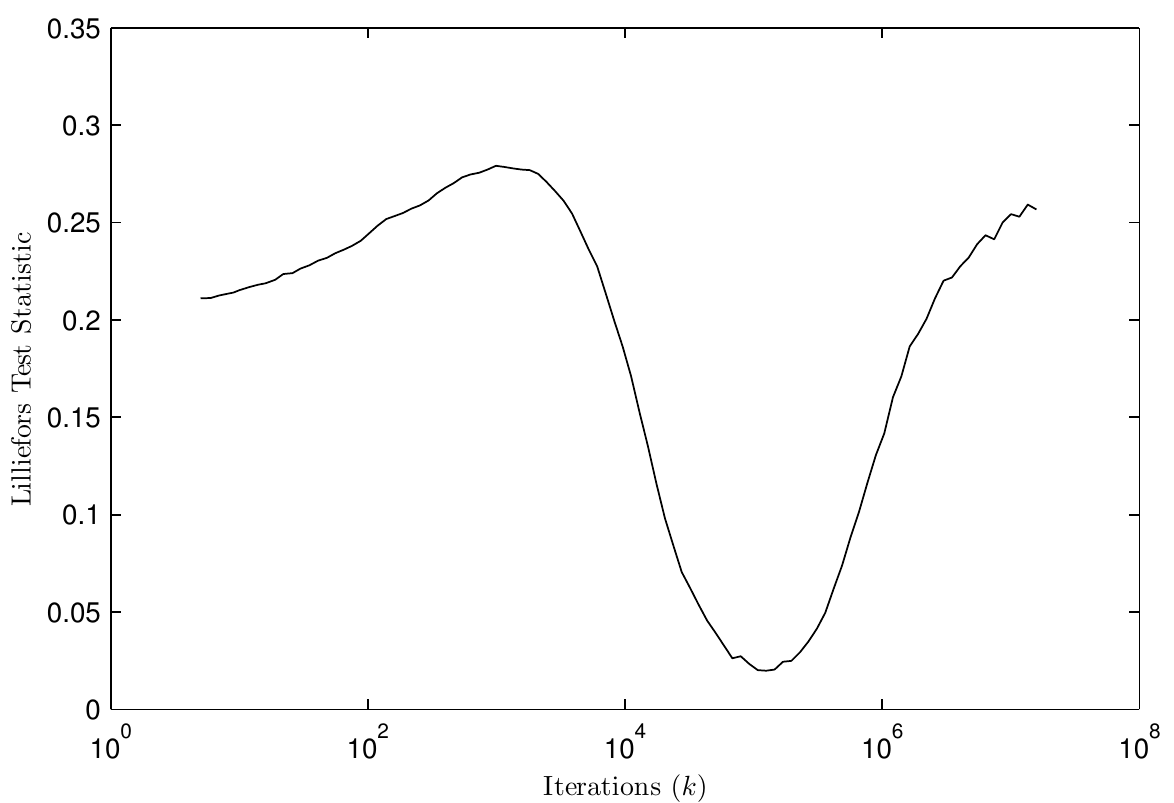}
\end{center}
\caption{The Lilliefors test statistic for the ensembles
$\{\theta(T^k(z_i))\}$
described in Algorithm \ref{algo-reduce}.
A lower test statistic indicates the distribution is closer to normal.
One iterate corresponds to a time step $\Delta t=\eng{2}{-6}$ of
the system \refsfwith.
}
\label{lillie-figure}
\end{figure}

Let us now comment on the choice of the flow time $t_k$ (or the number of iterations, $k$). The definition of the drift and diffusion coefficients requires $t_k \to 0$ \cite{Gardiner}. However, to capture time scales which go beyond the fast time scale, i.e., times which are sufficiently large to allow for equilibration on each fast fiber, $t_k$ needs to be sufficiently large. In the context of the particular scaling of the multiscale system \refsf\ we require $t_k \sim{\mathcal{O}}(1/\ep)$, which is long enough that the fast variables behave according to the invariant density on each fiber and small enough to allow for (nearly) constant slow variables.
For such times the diffusion term in the reduced dynamics is expected to be dominant; hence we chose $t_k$ such that the set of values $\theta(T^k(z_i))$ has
a distribution close to normal.
To test for normality we use the Lilliefors test for normality
\cite{lilliefors1987}.
In many cases, the Lilliefors test and the Kolmogorov-Smirnov test on which
it is based are used to test the hypothesis that a given finite set of
points was sampled from a normal distribution.  In our
situation, the projection of the dynamics only approximates a
normal distribution and it is not helpful to calculate the probability that
such a population came from a true normal distribution.
Instead, we consider
the test statistic used in the Lilliefors test, which is a positive number
computed from the sample, and use that this test statistic is closer to zero when the
sample is closer to normal. In Figure~\ref{lillie-figure} we show how $t_k$ can be chosen as the time where the Lilliefors test statistics has a minimum.\\

In the framework for slow-fast systems introduced in Section \ref{sec-sfabs}, the
drift-diffusion process \eqref{eqdd} now plays the role of the reduced
dynamical system $S$ on the reduced
space $\cX = \bbR$ and the map $\P:\cZ \to \cX$ is given by $\theta$.
By the assumption $\P \circ T^N \approx S \circ \P$ for large $N$, if $k$ is sufficiently
large, then for a point $z \in \theta \inv(\{v\})$, the value $\theta(T^k(z))$
gives an approximation of the result of evolving \eqref{eqdd} forward
by a time $t$ proportional to $k$ from the starting value $X_0=v$.
Thus, the computation
$
    v \leadsto z_i \leadsto \theta(T^k(z_i))
      $
in Step 5 of Algorithm \ref{algo-reduce} provides the ensemble
$\{X_{{t_k},i}\}$ described above.

\subsection{Examples} 

\subsubsection{The skew-product example} \label{sec-axis2} 

We now apply Algorithm \ref{algo-reduce} to the example system
\refsf\ from Section \ref{sec-axis}.
For this system, Figure \ref{lillie-figure} plots how
the test statistic varies in the number of iterates $k$ for a sample of points
on a fast fiber $F = \theta\inv(\{v\})$.
For $k = 0$, every point $z_i$ projects to the same
value in $\bbR$.
As $k$ increases, the distribution becomes more normal as the values
$\theta(T^k(z_i))$ spread out.  However, as $k$ becomes very large, the
distribution of the points starts to approach the invariant density of
the full system, and this density, when projected by $\theta$, will in general
not be normal.  From the graph, a value of $k=10^5$ iterates gives a
distribution which is close to normal.
That the test statistic is low here ($\sim 0.02$) gives empirical support to
the assumption that the slow dynamics can be well approximated by a
drift-diffusion process.

\begin{figure}[tp]
\begin{center}
\includegraphics{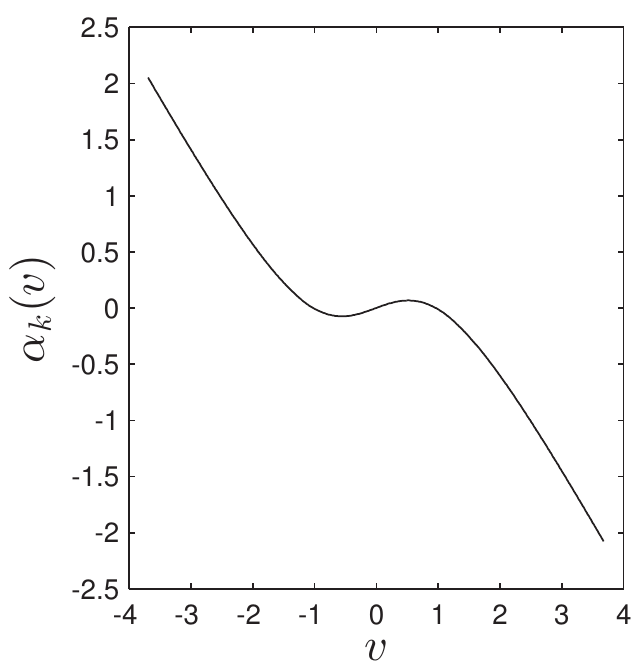}
\includegraphics{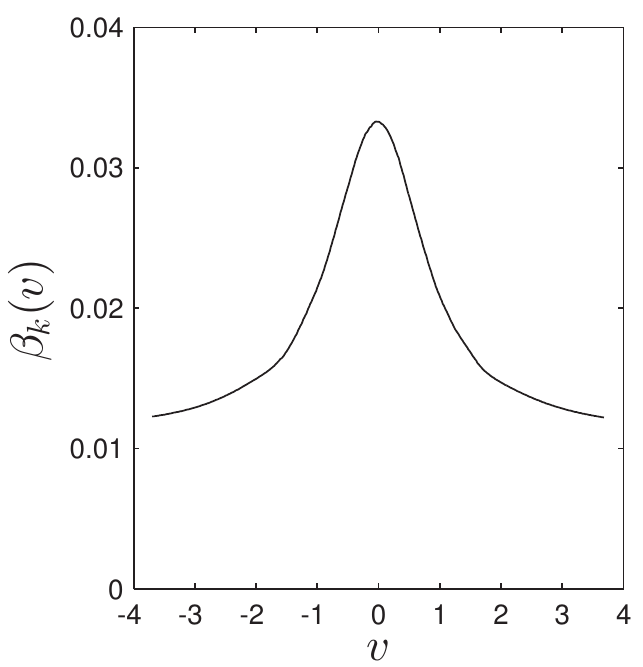}
\end{center}
\caption{The drift $\alpha_k$ and diffusion $\beta_k$ computed in Algorithm
\ref{algo-reduce} for the system \refsfwith.}
\label{alphabeta-figure}
\end{figure}

For simplicity, we assume that the same value $k$ can be used for all fibers,
as is the case for this example.
For Step 5 of Algorithm \ref{algo-reduce}, we used 500 fibers
with the value $v$
varying between $-4$ and $4$.
Each fiber was populated with $q = 10^4$ points.
The resulting functions $\alpha_k$ and $\beta_k$
are plotted in Figure \ref{alphabeta-figure}.

If $k$ iterations corresponds to a time interval $k \Delta t$, then $\alpha$
and $\beta$ in \eqref{eqdd} can be approximated by
${\alpha_k}/{k \Delta
t}$ and ${\beta_k}/{k \Delta t}$.
Note that while the diffusion function
$\beta$ would be constant for the homogenization limit of \refsf\ as $\ep \to
0$, the computed function $\beta_k$ in Figure \ref{alphabeta-figure} is
non-constant because it is computed using a non-zero time-step $k \Delta t$
which causes enhanced diffusion near the saddle point which separates the two
metastable states near $x = \pm 1$.\\

\begin{figure}[tp]
\begin{centering}
\includegraphics[scale=0.8]{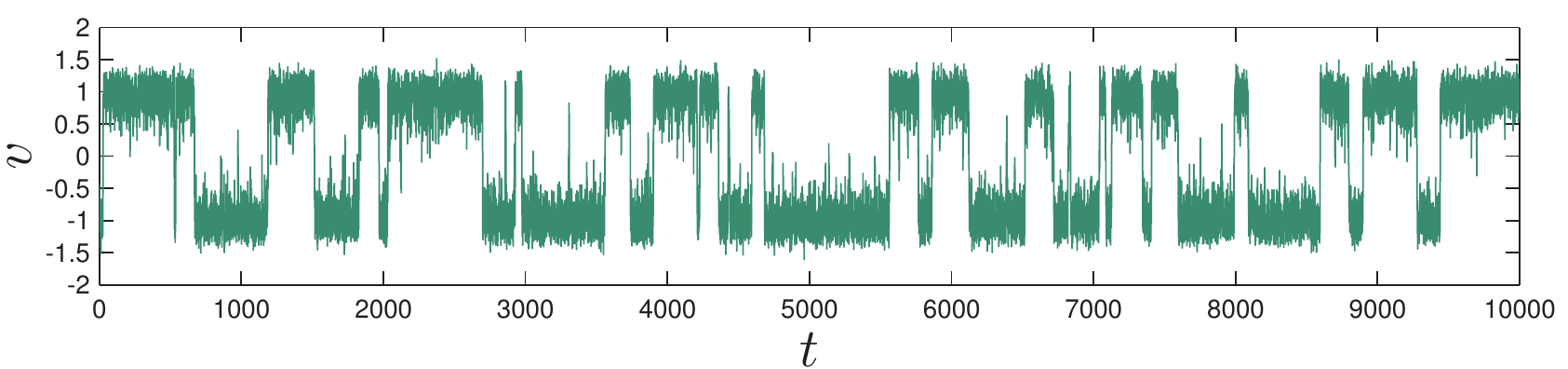}
\end{centering}
\caption{Simulation of an orbit of the reduced one-dimensional dynamics
given by Algorithm \ref{algo-reduce} for the system \refsfwith.
The orbit was calculated using \eqref{eqn-redstep} with $dt=0.02$.}
\label{fig-evolution}
\end{figure}

The drift-diffusion process \eqref{eqdd} can now be used to simulate the slow
dynamics using a time step which would be too large to simulate the full
dynamics
Here, we use a time step of $dt=0.02$ and calculate the orbit by
\begin{equation} \label{eqn-redstep}
    v_{n+1} = v_n + \alpha(v_n) \, dt + \sqrt{\beta(v_n)\, dt}\, r_n
\end{equation}
where the $r_n$ are independent pseudo-random values drawn from the standard
normal distribution.
Such a simulated orbit is plotted in Figure \ref{fig-evolution}.
This simulated orbit closely resembles the slow $x$ variable of the original
system (shown in Figure \ref{fig-fullevolution}).

\subsubsection{A distorted example} \label{sec-distort2} 
\begin{figure}[tp]
\begin{center}
\includegraphics{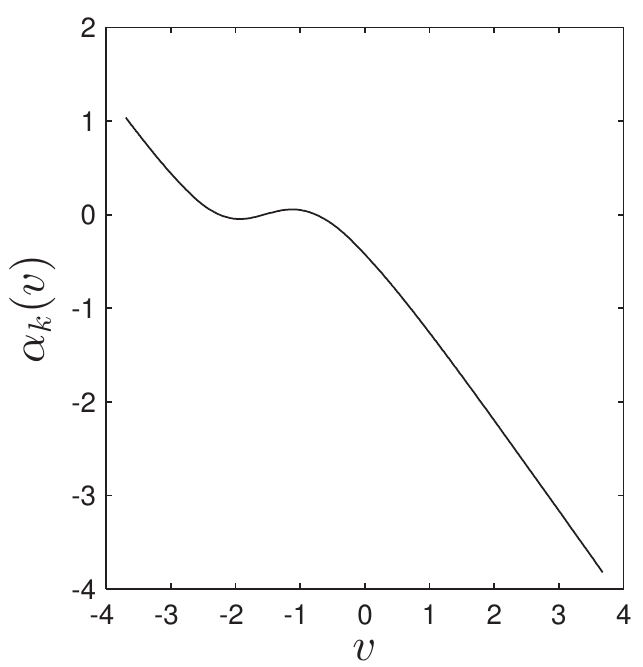}
\includegraphics{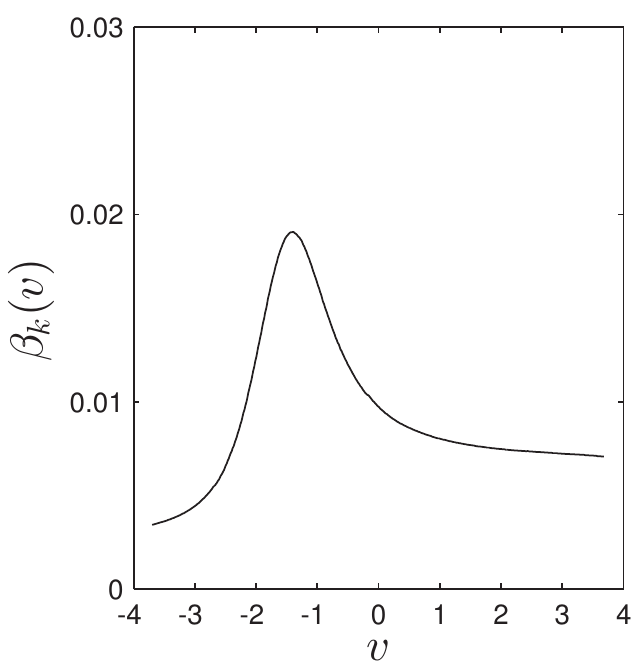}
\end{center}
\caption{The drift $\alpha_k$ and diffusion $\beta_k$ computed in Algorithm
\ref{algo-reduce} for system \refsfwith\ after
the distortion $g$ described in Section \ref{sec-distort}.}
\label{fig-robowalphabeta}
\end{figure}

\begin{figure}[tp]
\begin{centering}
\includegraphics[scale=0.8]{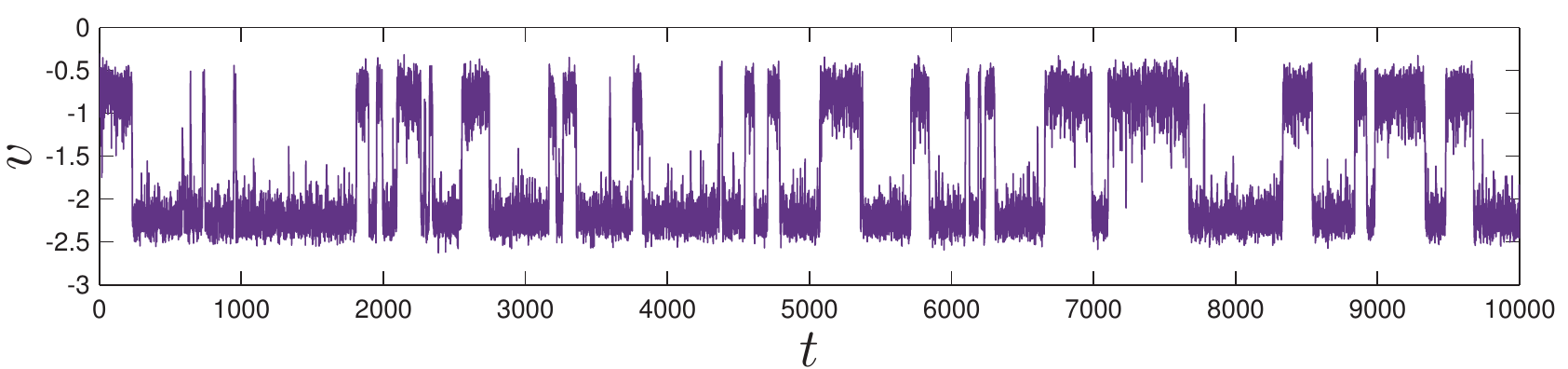}
\end{centering}
\caption{Simulation of an orbit of the reduced one-dimensional dynamics
given by Algorithm \ref{algo-reduce} for system \refsfwith\ after
the distortion $g$ described in Section \ref{sec-distort}.
The orbit was calculated using \eqref{eqn-redstep} with $dt=0.02$.}
\label{fig-robow-evolution}
\end{figure}

Finally, we use Algorithm \ref{algo-reduce} to simulate a slow orbit of the distorted system of Section \ref{sec-distort}.
The computed functions for the drift $\alpha_k$ and diffusion $\beta_k$ are
plotted in Figure \ref{fig-robowalphabeta} and a simulated orbit of the reduced
system is plotted in Figure \ref{fig-robow-evolution}. As in Section~\ref{subsec-dist1}, our algorithm is insensitive to non-trivial coordinate transformations and accurately distills the slow variables and their dynamics.

\subsection{Statistics} 

From a visual comparison of
Figures \ref{fig-fullevolution},
\ref{fig-evolution}, and \ref{fig-robow-evolution},
it is not entirely clear if the reduced dynamics of Algorithm \ref{algo-reduce}
is accurately capturing behaviour of the slow variable $x$ in \refsf.
To be more certain, we compute statistics related to these time series.

As argued in Section \ref{sec-finding}, the eigenfunction
$\phi$ of the Koopman
operator is approximately equal to $\psi \circ \P$ where $\psi$ is an
eigenfunction for the reduced dynamics and $\P:\cZ \to \cX$ is the mapping
down to the reduced space.  Algorithm \ref{algo-reduce} replaces $\phi$ with a
function $\theta$ which has similar level sets.  Therefore, $\theta$ should be
of the form $\theta = h \circ \P$ for some unknown $h$.  Because of this, when
comparing time series of the full dynamics with time series for the reduced
dynamics, one should only use statistics which are unchanged when $\P$ is
replaced with $h \circ \P$ where $h$ is a homeomorphism of the real line.

In our case, an obvious choice for such a statistic is the average switching
time between the two metastable states.
From a time series $\{v_i\}_{i=1}^M$, such as those plotted in
Figures \ref{fig-fullevolution},
\ref{fig-evolution}, and \ref{fig-robow-evolution},
compute percentiles $v^-$ and $v^+$ such that 40\% of the points in $\{v_i\}$
are below $v^-$ and 60\% are below $v^+$.  Then, declare each point $v_i$
to be in either the ``high'' state or the ``low'' state by the following
heuristic{:}
\begin{itemize}
    \item
    if $v_i  \ge  v^+$, then $v_i$ is in the high state,
    \item
    if $v_i  \le  v^-$, then $v_i$ is in the low state,
    \item
    otherwise, $v_i$ is in the same state as $v_{i-1}$.
\end{itemize}
This gives us a way to measure the switching times which is unaffected by
applying a homeomorphism $h:\bbR \to \bbR$ to the time series.

For the full dynamics \refsfwith, we computed a time series $\{x_i\}_{i=1}^M$
of the $x$-variable with $M = \eng{2}{6}$ where each step $x_i$ to
$x_{i+1}$ represents a
change of time of $\tau = 0.2$.  This involved $\eng{2}{11}$ individual steps
of the Euler-Maruyama method, using $\Delta t = \eng{2}{6}$.
For the reduced dynamics, a similar time series
$\{v_i\}_{i=1}^M$ with the same $M=\eng{2}{6}$ and $\tau=0.2$ was obtained.
Each step $v_i$ to $v_{i+1}$ was computed by ten applications of
\eqref{eqn-redstep}, each with a time step of $dt = 0.02$.

For the full time series, there were 2222 switches between the
slow metastable
states, with a mean switching time of $\tau_0 = 180.0$.
For the reduced dynamics of Section~\ref{sec-axis2}, there were 2280 switches
with a mean time of $\tau_0 = 175.5$.
In the distorted case in Section \ref{sec-distort2}, there were 2178 switches
with a mean time of $\tau_0 = 183.7$.

\begin{figure}[tp]
\begin{centering}
\includegraphics{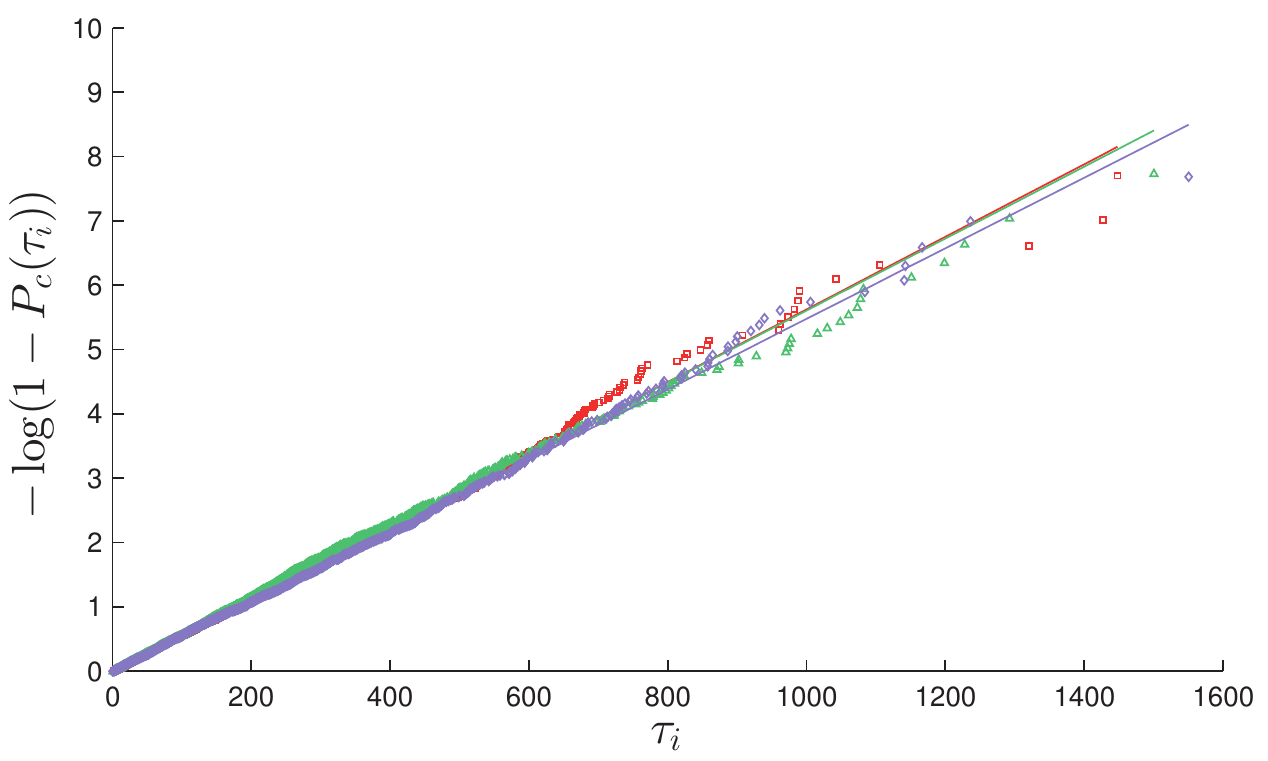}
\end{centering}
\caption{Log plot of the normalized histogram of switching times $\tau_i$.
This is computed for time series of the $x$ variable of \refsf\ as in
Section~\ref{sec-axis} (red squares),
for the reduced dynamics in Section~\ref{sec-axis2} (green triangles)
and for the reduced dynamics of the distorted system in Section~\ref{sec-distort2}
(blue diamonds).}
\label{fig-gaps}
\end{figure}
In each case, we expect the switching times $\tau_i$ to obey a Poisson process,
where the cumulative probability function is given by the formula
\[
    P_c(\tau_i) = 1 - \exp(\frac{\tau_i}{\tau_0}),
\]
where $\tau_0$ is the mean switching time.
To see if this the case, the computed
switching times $\tau_i$ are plotted
against $-\log(1-P_c(\tau_i))$, as shown in Figure \ref{fig-gaps}.
For each of the three time series, the resulting plot is linear,
and a least squares fit gives a computed value for $\tau_0=177.3$ for the full dynamics $\{x_i\}$ of Section~\ref{sec-axis},
$\tau_0=178.8$ for the reduced dynamics of Section~\ref{sec-axis2}, and
$\tau_0=182.2$ for the reduced dynamics of the distorted system in
Section~\ref{sec-distort2}.
This suggests that Algorithm 2 generates dynamics which closely resemble the
relevant dynamics of the full system.


\section{Computing the eigenfunctions} \label{sec-computing} 

The techniques in this paper use spectral properties of the transfer
operator, and its adjoint, the Koopman operator. Therefore, we
must approximate these operators numerically. There are a number of
approaches that could be used, and specialized methods other than Ulam, may be significantly faster in certain
settings, depending on the smoothness of the system.
These include numerical methods based on the infinitesimal generator
of Ulam's method \cite{fjk2013},
finite-element methods \cite{BabuskaOsborn89},
and Galerkin methods using Fourier or Chebyshev bases \cite{Boyd01}.
See also Chapter 5 of \cite{KoltaiThesis}.
For simplicity, we employ the standard Ulam's method
\cite{ulam-collection, li1976} for all of the examples in this paper. Eigenpairs and the level sets were computed using the {\tt eigs} and {\tt
contour} routines of MATLAB \cite{Matlab}.\\

Recall that for a deterministic dynamical system $T:\cZ \to \cZ$, the
transfer (or Perron-Frobenius) operator $\cL$ on $L^2(\cZ)$ is defined by
\begin{equation} \label{eqn-perfrob}
    \int g \cdot \cL f \,dm  = \int (g \circ T) \cdot f \,dm
\end{equation}
for all $f,g$ in $L^2(\cZ)$.
To solve numerically, Ulam's method divides the domain $\cZ$ into a finite
number of boxes $B_1, \ldots, B_n$ and reduces the problem to
\begin{equation} \label{eqn-galerkulam}
    \int \bbOne_{B_i} \cdot \cL \bbOne_{B_j} \,dm = \int (\bbOne_{B_i} \circ T) \cdot \bbOne_{B_j} \,dm
\end{equation}
for $i,j = 1, \ldots, n$ where $\bbOne_{B_i}$ is the indicator
function for the set $B_i$.  The resulting operator, now defined on the finite
dimensional space spanned by the functions $\bbOne_{B_i}$,
is equivalent to an $n
\times n$ matrix defined by
\[
    P_{ij} = \frac {m(T \inv(B_i) \cap B_j)} {m(B_i)}.
\]
To estimate $P_{ij}$ numerically, populate each $B_j$ with a large number
$N$ of randomly sampled test points $z_{j,k}$ for $k=1,\ldots,N$ and
approximate
\begin{equation} \label{eqPijapprox}
    P_{ij} \approx \frac {1}{N} \# \{k : T(z_{j,k}) \in B_i \}.
\end{equation}
In the case of a random dynamical system, the dynamics is defined by a
map $T:\cZ \times \Omega \to \cZ$ where $\Omega$ is a probability space, sometimes
called the ``noise space'' of the system.  Each $\omega \in \Omega$ defines
a function $T_\omega = T(\cdot, \omega)$ and one can define a transfer operator
$\cL_\omega$ for each $T_\omega$.
We then consider the annealed (or integrated) operator $\cL:L^2(\cZ) \to
L^2(\cZ)$ defined by
\[
    \cL(f) = \int_\Omega \cL_\omega(f) d \omega.
\]
Similarly, the annealed Koopman operator is
\[
    (\K f)(z) = \int_\Omega f(T(z,\omega)) d \omega
\]
which for a point $z$, gives the expected value of $f(T(z))$.
For a random dynamical system,
Ulam's method
%
produces a matrix which
can be approximated numerically by a slightly modified version of
\eqref{eqPijapprox} with
\begin{equation} \label{eqPijapproxRandom}
    P_{ij} \approx \frac {1}{N} \# \{k : T_{\omega_{j,k}}(z_{j,k}) \in B_i \}\, ,
\end{equation}
where for each $z_{j,k} \in B_j$
an element $\omega_{j,k} \in \Omega$ of the noise space
is chosen randomly and independently.

The examples in this paper are defined by stochastic differential equations.
In this case, the noise space $\Omega$ is the space of continuous paths
$C^0(\bbR, \bbR^d)$ where the probability measure is given by the Wiener process
representing Brownian motion.  Once a path $W \in \Omega = C^0(\bbR,\bbR^d)$ is
chosen, an SDE of the form
\[
    dx_t = \mu(x) dt + \sigma dW_t
\]
at starting value $x(0)=x_0$ is solved by a path satisfying the integral
\[
    x(t) = x_0 + \int_{0}^{t} \mu(x(s))\, ds + \sigma W(t)
\]
at all times $t > 0$.  For a specified flow time $\tau > 0$, this gives a
deterministic map $T_W(x_0)  = x(\tau)$ for each $W \in \Omega$ and so defines
a random dynamical system in the form $T:\cZ \times \Omega \to \cZ$.
For each point $z_{j,k}$ in \eqref{eqPijapproxRandom}, the image
$T_{\omega_{j,k}}(z_{j,k})$ can be computed, for example, by the Euler-Maruyama method \cite{kloeden1995},
and so the annealed operators for these
types of systems can be approximated numerically.

%
%

For these continuous-time dynamical systems,
one must choose some interval of time $\tau > 0$ in order to apply Ulam's
method to the map $T=\varphi_\tau$ given by the flow $\{\varphi_t\}_{t \in \bbR}$.
Picking too small a flow time can lead to a dynamical system appearing to
have more diffusion than is actually present, a problem closely related to
so-called ``numerical diffusion'' which appears in finite-difference methods
\cite{fjk2013}.
A common practice is to choose $\tau$ large enough, so that for a
significant number of sample points $x_{i,k} \in B_i$, the image
$T(x_{i,k})=\varphi_\tau(x_{i,k})$
is no longer inside the box $B_i$ \cite{froypad2009}.
In other words, so that $P_{ii}$ is not too close to one.

This works well for many systems, but it is not a sufficient criterion for a
system with multiscale behaviour.
For these systems, choosing a flow time in such manner
only guarantees that Ulam's method captures the fastest dynamics, and there
may be significant numerical diffusion present in the direction of the
slow dynamics, meaning that the computed invariant density and other
eigenfunctions may have significant errors.

The effect is more pronounced when the slow and fast directions do not align
with the grid of boxes.  Here,
Ulam's method cannot faithfully capture motion happening purely in the fast
direction.  In some sense, the Galerkin projection ``smears'' the point
$T(x_{i,k})$ over the box $B_j$ in which it resides and this smearing happens
in both slow and fast directions.
As an example, we use Ulam's method for the system \refsf\ after distortion as defined in Section
\ref{sec-distort} with an extremely short time-step of $\tau = \eng{2}{-5}$.
\begin{figure}[t]
\includegraphics{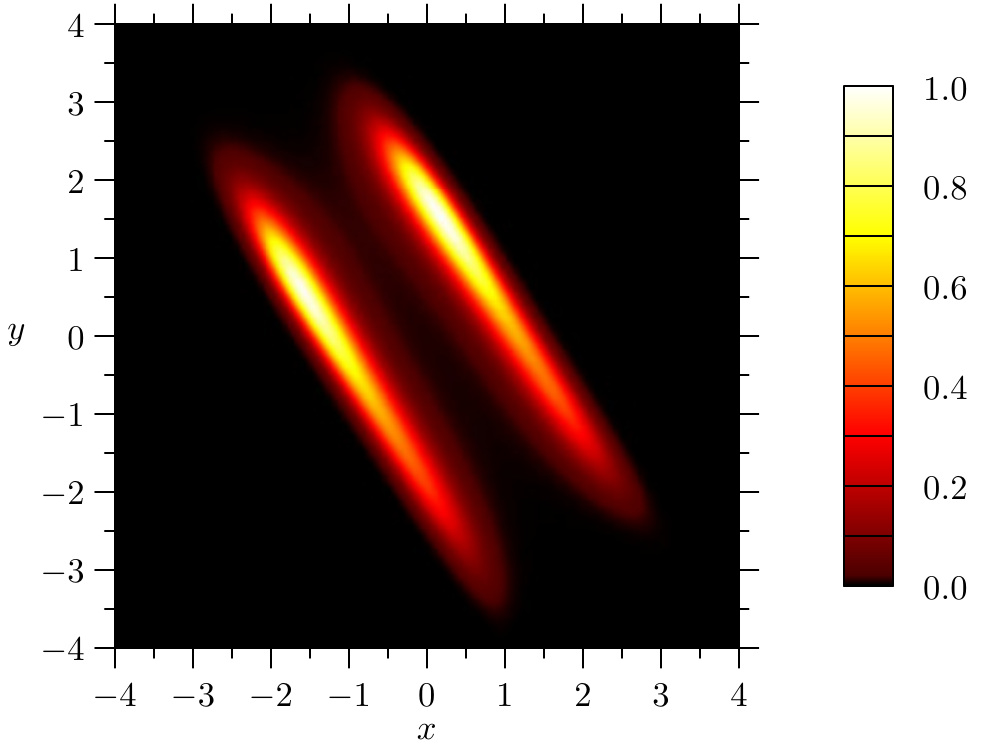}
\caption{An incorrect density calculated using Ulam's method with too short a
flow time for system \refsfwith\ after
the distortion $g$ defined in Section~\ref{sec-distort}.
The flow time used to construct the Ulam
matrix is $\tau=\eng{2}{-5}$. Compare this to the computed density in
Figure \ref{fig-dist-density}.}
\label{fig-numdiff}
\end{figure}
In the resulting Ulam matrix, the maximum value $P_{ii}$ on the diagonal is
less than $0.11$.
However, the computed density, plotted in Figure \ref{fig-numdiff},
shows two thick spurious bands, which would be more consistent with a system with
a large amount of diffusion in the slow direction (cf. the invariant density in Figure~\ref{fig-dist-density} computed with significantly larger flow time $\tau=0.04$).

One way to avoid this is to increase the flow time so that a typical point
moves in the slow direction by an amount comparable to the size of the boxes
$B_i$.  However, we are assuming that we do not know in advance the directions
of the slow and fast dynamics.
To handle this, we compute unit eigenvectors using Ulam's method for a
sequence of increasing flow times, and stop when the result stabilizes.
For instance, in the axis-aligned example system \refsf\ in Section
\ref{sec-axis}, the invariant density
$\rho_\tau \in \spann\{\chi_1, \cdots, \chi_n\} \subset L^1(\bbR^2)$
was computed for flow times $\tau = 0.01, 0.02,$ and $0.04$.
Using the $L^1$ norm $\|f\| = \int_{\bbR^2} |f(x)| dx$, we compute
$\|\rho_{0.01} - \rho_{0.02}\| \cong 0.08$ and
$\|\rho_{0.02} - \rho_{0.04}\| \cong 0.03$.
More important than the density is the second leading eigenfunction of the
Koopman operator,
i.e., the eigenfunction shown in Figure \ref{evec-figure},
since this function is used to define the fiber dynamics.
The corresponding $L^1$ differences for this computed eigenfunction are
$\|\phi_{0.01} - \phi_{0.02}\| \cong 0.04$ and
$\|\phi_{0.02} - \phi_{0.04}\| \cong 0.01$.
We view this difference as small enough that we can safely use $\tau = 0.04$
in the computation of the eigenfunctions.

Because we are using the Euler-Maruyama method with a time-step of
$\eng{2}{-6}$, this value of $\tau$ corresponds to $\eng{2}{4}$ iterations.
However, this value of $\tau$ is still much smaller than the average time an
orbit takes to switch between the metastable states based at $x = -1$ and $x =
+1$.  That switching happens on the order of hundreds of units of time
(as shown in Figure \ref{fig-fullevolution}).

As a further test, we also consider the computed eigenvalues.  As argued in
Section \ref{sec-ops}, the leading eigenvalues of the full dynamics
correspond to decay rates associated to the slow dynamics.
Therefore, if many of the computed eigenvalues are very close to one, it
suggests that the flow time may be too short.
For $\tau = 0.04$ in the example system, the leading ten computed eigenvalues
(as listed in Table \ref{table-sf}) range from 1.000 to 0.6969
which is far enough from 1 to tell us that
advection in the slow direction is being captured in the Ulam matrix.\\

In addition to choosing the flow time, one must also set the number of sample
points $N$ in \eqref{eqPijapprox} or \eqref{eqPijapproxRandom}.
Because of the relatively long flow time, starting
from one box $P_{ij}$ the evolved points $T(x_{i,k})$ may lie in a large
number of boxes $B_j$, and $N$ must be chosen large enough that each $P_{ij}$
is computed to a reasonable accuracy.
In all our examples, we use $N = 10^4$ and a 200 by 200 grid
on the region $[-4, 4] \times [-4, 4] \subset \bbR^2$.  This is a much
larger region than where orbits typically lie.  However, we compute everything
in this larger region so that when fiber dynamics are computed,
the orbit will stay in the domain with very high probability.

Once a level set is computed, the same heuristics using $L^1$-distance of the
computed eigenfunctions and the distance of the eigenvalues from $1$ may be
used to determine a reasonable number of iterations to apply for the fiber
dynamics.
These showed that 20 iterations for a total flow time of $20 \times
\eng{2}{-6} = \eng{4}{-5}$ was sufficient for the fiber dynamics of the
example system \refsf\ in Section \ref{sec-axis}.

An eigenfunction $\phi$ computed using Ulam's method will be given
as a linear combination of the indicator
functions $\bbOne_{B_i}$.
As such functions are piecewise constant, the gradient $\nabla \phi$ will
therefore be zero anywhere it is defined.
This gradient is not useful in defining a gradient flow $\pi$ as used in
Step 4 of Algorithm \ref{algo-compare}.
Instead, when computing $\nabla \phi$ for use in the algorithm
we replace $\phi$ by the bilinear
interpolation of values it takes on the centers of each box $B_i$.
Further, as solving for the gradient flow is computationally expensive,
we employ the following shortcut.
To compute the map $\pi$ at a point $z$ close to the fast fiber,
evaluate the gradient at $z$, then take $\pi(z)$ to be the intersection of
the fiber $\phi \inv(\{v\})$ with the line given by $z + \bbR \cdot \nabla \phi(z)$.
Since in computing the fiber dynamics the function $\pi$ is only evaluated
at points very close to the level set, the gradient is nearly constant on the
short path from $z$ to $\pi(z)$ and so the approximation is valid.\\

Alternatively, at least in principle, one could have used the eigenfunctions of the transfer operator $\cL$ instead of those of the Koopman operator $\K$. If $\mu$ is an invariant density for
the dynamical system $T$, then \eqref{eqn-perfrob} with $m=\mu$
defines a transfer operator $\cL_\mu$ and the invariance of the measure implies
that a constant function $f \equiv c \in \bbR$ satisfies $\cL_\mu(f) = f$.
Numerical experiments suggest that, similar to the Koopman operator, leading
eigenfunctions of $\cL_\mu$ are nearly constant in the fast direction of a
multiscale system.
The problem with this approach is that one must numerically approximate the
invariant measure $\mu$,
say by a density function, $\rho$, and use this to compute
eigenfunctions of $\cL_\mu$.  Numerically, this is equivalent to calculating
the pointwise quotient $\phi/\rho$ where $\phi$ is an eigenvalue of the
transfer operator $\cL$ defined with respect to Lebesgue measure $m$.
For systems which involve regions of low density with $\rho$ close to zero, this is intractable to do numerically.

\section{Discussion and outlook} 
\label{sec-discussion}
We have used the transfer operator and the Koopman operator to disentangle the multiple slow and fast scales and provide an effective set of slow variables which can then be efficiently numerically integrated in time. We constructed a projection onto a slow subspace using the eigenfunctions of the Koopman operator. We then devised an algorithm to estimate the drift and diffusion terms of an assumed stochastic differential equation describing the reduced slow dynamics. This constitutes a huge computational advantage as we were able to use time steps of the order of $10^4$ larger to propagate forward the slow variables when compared to the full multiscale system.\\


There are similarities of our approach with isocommittor surfaces in transition path theory (see for example \cite{EVandenEijnden04} for a review) and with reaction coordinates for constrained molecular dynamics (see for example \cite{HartmannSchuette05,HartmannSchuette07,Lelievre,LelievreEtAl12}). Isocommittor surfaces separate metastable states and these methods aim at calculating the associated free energy in order to determine statistics such as most likely transitions and transition rates between the meta-stable states. It was, however, noted in \cite{E} that reaction coordinates and isocommittor functions are not necessarily slow variables. In our method the slow and fast variables are first identified, and then subsequently the reduced dynamics is provided for the slow variables, whereas reaction coordinates and isocommittor surfaces are constructed to obtain statistics of transitions between predefined sets such as meta-stable states.\\

In the current paper,
the algorithms were described and tested in a low-dimensional setting
with one-dimensional slow dynamics and two time scales.
As an outlook for future research,
we now discuss
how one might extend these techniques to more complicated situations.

\smallskip

\noindent
{\bf Extension to systems with multi-dimensional slow subspace:} In this paper, we considered multiscale examples with one-dimensional slow
dynamics, so that the mapping $\cP:\cZ \to \cX$ could be approximated by
a single eigenfunction $\phi:\cZ \to \bbR$ of the Koopman operator.
We outline now a possible extension
of our computational method to the case of a
higher-dimensional slow subspace.  Details are planned for further research.
For a system with higher-dimensional slow dynamics, it will be necessary
to approximate $\cP$ by a product of eigenfunctions
\[
    \cP = \phi_1 \times \phi_2 \times \cdots \times \phi_m : \cZ \to \bbR^m.
\]
Then, a value $v \in \bbR^m$ defines a fiber $F = \cP \inv(\{v\}) \subset \cZ$.
The fiber dynamics
can be defined as $\hat T:F \to F$ where $\hat T = \pi \circ T$ and
$\pi:U \to F$ is a mapping from a neighbourhood $U$ down to the
fiber itself.
As explained in Section \ref{sec-computing},
for a single eigenfunction $\phi$
the numerical computation of $\pi(z)$ was given by the
intersection of $F$ with the line through $z$ tangent to $\nabla \phi(z)$.
A possible way to extend this to the general setting
is to define
$\pi(z)$ as the intersection of $F$ with the $m$-dimensional affine subspace
through $z$
spanned by the vectors $\nabla \phi_i(z)$ for $i = 1,\ldots,m$.
Then, as in Section \ref{sec-axis},
one could test the validity of the fiber dynamics by
comparing the distance $\|\hat T(z) - T(z)\|$ to the distance
$\|T(z) - z\|$ for an ensemble of points on $F$.

If the dimension of the slow dynamics is not known beforehand,
an automated method to find the dimension would be to test
products of an increasing number of eigenfunctions until a product is found
for which $\hat T$ gives a reasonable approximation of fast dynamics.
Using this definition of $\hat T$, Algorithm \ref{algo-compare}
could be implemented much as before.

Algorithm \ref{algo-reduce} might also be generalized to use a product of
eigenvalues.  In such a case,
the
$\alpha_k(v)$ would be computed as the
$m$-dimensional mean of an
ensemble
of the form
$\{\theta(T^k(z_i)) - v\}_{i=1}^q$
using a product
$
    \theta = \theta_1 \times \cdots \times \theta_m
$
and $\beta_k(v)$ would be computed as the $m \times m$ covariance matrix.

The main difficulty which would arise
in implementing these extensions of the algorithms to
high-dimensional systems
is accurate and efficient numerical approximation of the transfer operator.
Currently, techniques based on Ulam's method are mainly restricted to
low-dimensional systems.

\smallskip

\noindent {\bf Extension to systems with more than two time scales:} In this
paper, we considered dynamical systems with one fast and one slow time scale.
We now suggest an iterative procedure using Algorithms~\ref{algo-compare} and
\ref{algo-reduce} (or their higher-dimensional variants) which distills the
dynamics on several time scales. For ease of exposition, we present this
procedure for a system $T$ with three time scales (fast, moderate and slow).
Applying Algorithm~\ref{algo-compare} to such a system produces
fiber dynamics $\hat T:F \to F$ which incorporates the dynamics of both the
fast and moderate time scales.
The slowest time scales can be obtained by applying Algorithm \ref{algo-reduce}
to $T$.
To further isolate the moderate and fast dynamics,
apply Algorithms \ref{algo-compare} and \ref{algo-reduce} now to the fiber
dynamics $\hat T$ in place of the full system $T$.

\smallskip

\noindent
{\bf Solving for the slow variables without having to integrate forward in
time the fast variables:} We defined the projection to the slow subspace via
the eigenfunctions of the Koopman operator. This requires the integration in
time of the full dynamics, including the fast variables, needed for Ulam's
method. The transfer operator $\cL$ and the infinitesimal
generator $\A$ share the same eigenfunctions
and this has been exploited to find these eigenfunctions
numerically for certain systems without trajectory integration
\cite{fjk2013,KoltaiThesis}.
Thus, it may be possible to reformulate the steps in Algorithm 1 and 2
based purely on analysing the vector field itself without numerical
integration, yielding a significant speed up in computation.


\section*{Acknowledgments}
The illustrations in Figures \ref{fig-Proj} and \ref{fig-Reduce} were produced
using the Asymptote vector-graphics language ({\tt http://asymptote.sf.net}).
We acknowledge the support by the Australian Research Council under
grant DP$120104514$.

\bibliographystyle{siam}
\bibliography{multi}

\end{document}